\documentclass{amsart}
\usepackage{color}
\usepackage{amsmath}
\usepackage{amsthm}
\usepackage{amssymb}
\usepackage{graphicx}
\usepackage{amsfonts}
\usepackage{mathrsfs}

\theoremstyle{plain}
\newtheorem{theorem}{Theorem}[section]
\newtheorem{lemma}[theorem]{Lemma}

\theoremstyle{definition}

\newtheorem{remark}[theorem]{Remark}
\newtheorem*{remark*}{Remark}

\def\dbar{\overline{\partial}}
\def\pmtwo#1#2#3#4{\left( \begin{array}{cc}#1&#2\\#3&#4\end{array}\right)}

\begin{document}

\title{Spectral approach to D-bar problems}

\author[C.~Klein]{Christian Klein}
\address[C.~Klein]
{Institut de Math\'ematiques de Bourgogne
9 avenue Alain Savary, BP 47870, 21078 Dijon Cedex}
\email{christian.klein@u-bourgogne.fr}

\author[K.~McLaughlin]{Ken McLaughlin}

\address[K.~McLaughlin]
{Department of Mathematics, The University of Arizona 617 N. Santa Rita Ave. P.O. Box 210089 Tucson, AZ 85721-0089 USA}
\email{mcl@math.arizona.edu}

\begin{abstract}
We present the first numerical approach to D-bar problems having 
spectral convergence for real analytic rapidly decreasing potentials. The 
proposed method starts from a formulation of the problem in terms of an 
integral equation which is solved with Fourier techniques. The 
singular integrand is regularized analytically. The resulting 
integral equation is  approximated via a discrete system which is 
solved with Krylov methods. As an example, the D-bar problem for the 
Davey-Stewartson II equations is solved.  The result is used to test 
direct 
numerical solutions of the PDE. 
\end{abstract}

\date{\today}

\subjclass[2000]{}

\thanks{}
\maketitle

\section{Introduction}
It is the purpose of this paper to present an efficient numerical approach to 
solve D-Bar equations for a complex valued unknown $\psi = 
\psi(x,y)$, i.e.,
\begin{eqnarray}
\label{eq:001}
\dbar \psi = \frac{1}{2}Q(x,y) \overline{\psi}.
\end{eqnarray}
In (\ref{eq:001}), one has $(x,y) \in 
\mathbb{R}^{2}$, and we often use $z = x + i y$, as well as  $\dbar = 
\frac{1}{2} \left( \frac{\partial}{\partial x} + i  
\frac{\partial}{\partial y} \right)$. The function $Q$ is a given 
``potential'' which we assume to be a function in the Schwartz class 
$\mathcal{S}(\mathbb{R}^{2})$ 
of rapidly decreasing smooth functions in the plane.

This equation appears in a number of different areas of analysis and 
applied mathematics, most prominently in Electrical Impedance Tomography  which is of enormous importance in medical 
imaging, in the solution of completely integrable $2+1$ dimensional 
partial differential equations (PDEs), as well as in two-dimensional 
orthogonal polynomials and random matrix theory (see below for 
details and references). 


In this paper we will primarily be interested in  
\emph{Complex Geometrical Optics} (CGO) solutions of  (\ref{eq:001}) 
defined by the  normalization at infinity,
\begin{eqnarray}
\psi e^{- kz}  = 1 + \mathcal{O} \left(  \frac{1}{|z| } \right) \ ,
\end{eqnarray}
where $k \in \mathbb{C}$ is an auxiliary parameter.

In order to make use of the Fourier transform, we express $\psi$ in terms of $m$ which tends to $0$ as $|z| \to \infty$:
\begin{eqnarray}
\label{eq:038}
m = \psi e^{- kz } - 1 \ .
\end{eqnarray}
This new quantity now solves
\begin{eqnarray}
\label{eq:039}
\dbar m = \frac{1}{2} e^{\overline{kz}-kz} Q (\overline{m} + 1) \ .
\end{eqnarray}

Taking the Fourier transform defined via
\begin{eqnarray}&& \ \ 
\mathcal{F}(f) = 
\hat{f}(\xi) = \frac{1}{2 \pi} \int_{\mathbb{R}^{2}} f(x,y) e^{- i \xi_{1}x - i \xi_{2} y } dx dy \ , \ \ \ \ \mathcal{F}^{-1}(\hat{f}) = f(x,y) = \frac{1}{ 2 \pi} \int_{\mathbb{R}^{2}} \hat{f}(\xi_{1},\xi_{2}) e^{ i \xi_{1} x  + i \xi_{2} y} d \xi \ ,
\end{eqnarray}
we find that (formally) $\hat{m}$ solves 
\begin{eqnarray}
\label{eq:08}
 i  \xi \hat{m} =  \mathcal{F}\left\{  Q e^{ \overline{kz}-kz} \left( \overline{m} \right) \right\} 
 + \mathcal{F}\left\{  Q e^{ \overline{kz}-kz}  \right\}  \ .
\end{eqnarray}
It turns out to be convenient for the numerical approach to write the equation in terms of 
\begin{eqnarray}
\label{eq:042}
S(\xi ) = \xi \ \hat{m}(\xi) \ , 
\end{eqnarray}
since $\hat{m}$ is unbounded for $|\xi|\to0$, whereas $S$ is not. Thus the fundamental equation of study is
\begin{eqnarray}
\label{eq:MainEQa}
S(\xi) =  - i \mathcal{F}\left\{  Q e^{ \overline{kz}-kz} \left(\overline{ \mathcal{F}^{-1} \left( \frac{1}{\xi} S(\xi) \right) }  \right) \right\} 
 -i  \mathcal{F}\left\{  Q e^{ \overline{kz}-kz}  \right\} \ .
\end{eqnarray}

A brief gaze at (\ref{eq:MainEQa})  reveals the 
primary obstacle for a numerical approach:  the division by $\xi$, although a weak singularity in two dimensions, needs to be addressed.  

The approach we have developed to overcome this obstacle is to 
re-write the equation in a manner which peels off the singular 
behavior.  The offending term appearing in (\ref{eq:MainEQa}) is re-written as follows:
\begin{eqnarray*}
\mathcal{F}^{-1} \left( \frac{1}{\xi} S(\xi) \right) = \mathcal{F}^{-1} \left( \frac{1}{\xi} \left( S(\xi)  - G(\xi) \right) \right)  + \mathcal{F}^{-1} \left( \frac{1}{\xi} G(\xi) \right) \ , 
\end{eqnarray*}
and the function $G$ is chosen so that the following conditions are met:
\begin{itemize}
\item The first term above is smooth at $\xi=0$ (to machine precision) and vanishing (again to machine precision) at the boundary of the computational domain.
\item The second piece, containing all the singular terms, is evaluated via {\it explicit 
analytical calculation.} 
\end{itemize}
Below we choose the function $G$ to essentially contain the "anti-holomorphic part" of the function $S$, or more precisely a truncated Taylor series representing the anti-holomorphic part of $S$, multiplied by a Gaussian.

For the numerical evaluation, the Schwartz functions are treated as essentially periodic (the computational domain is chosen large enough that the studied 
functions and their derivatives vanish to within numerical precision at the 
boundaries of the domain). The Fourier series is approximated 
numerically via a discrete Fourier series, computed via a Fast Fourier 
Transform (FFT). Thus the integral equation (\ref{eq:MainEQa}) is 
approximated via  a system of linear equations for a matrix $S$ (denoted with 
the same symbol in an abuse of notation) of finite dimension. 

\vskip 0.3in
In many applications, including the Davey-Stewartson equation, one often needs the solution for many values of the auxiliary parameter $k$.  The approach briefly outlined above works for small values of $|k|$, and the ensuing  system can be solved iteratively via GMRES \cite{gmres}.   

We show, however, that for large values of $k$, this first approach does not work adequately, because the support of the desired solution is no longer distant from the boundary of the computational domain (for the Fourier transform variable).  We develop a re-formulation of the (Fourier transformed) $\dbar$ equation based on the support properties of the solution for large $k$, which yields a once iterated variant of (\ref{eq:MainEQa}) and handles disparate support issues via careful shifting within the computational domain.  Remarkably, this new equation is equivalent (up to a Fourier transform) of the iterated system of equations which Perry studies in (\cite{Perry2012}), and in fact can be effectively used for {\it all values of $k$}.  

We study the convergence of the method for concrete examples and show that the numerical approach has \emph{spectral convergence}, i.e., that the numerical 
error decreases exponentially for smooth functions. The approach is 
then applied to the D-bar problem of the Davey-Stewartson (DS) II 
equation, and the resulting solution is compared to a direct 
numerical solution of DS II as in \cite{KR2011}.

The paper is organized as follows: In section 2 we briefly summarize some of the main applications of the D-bar problem, and 
previous numerical approaches to its solution. In section 3 we present 
a short review on the existence theory for equation 
(\ref{eq:MainEQa}). In section 4 we discuss  
the numerical implementation of a regularization approach for 
$|\xi|\to0$ of the integrand in (\ref{eq:MainEQa}). 
In section 5 we derive a once iterated version of 
equation (\ref{eq:MainEQa}) which is shown to be convenient for the 
construction of CGO solutions to D-bar problems, and use the approach of section 4 (along with careful shifting) for this equation. 
This approach  is applied to the D-bar problem for the DS II system as an 
example in section 6. The obtained solution is compared to a DS II 
solution for the same initial data via a direct numerical integration 
of DS II. We add some concluding remarks in section 7. 

\section{Applications of D-bar problems and numerical approaches}
In this section we summarize the main applications of D-bar problems, 
and review previous numerical approaches for its solution. 

\subsection{Application 1:  The inverse conductivity problem}

Briefly, one is given static electric measurements on the boundary of a domain, and the problem is to recover the conductivity in the domain.  This is used, for example, in medical Electrical Impedance Tomography.  The practical implementation involves tackling several different mathematical problems, one of which is precisely the $\dbar$ problem which is the topic of this paper.  The mathematical formulation centers around the so-called conductivity equation,
\begin{eqnarray}
\label{eq:Conduct1}
\nabla \cdot \left(  
\gamma(x,y) \nabla u(x,y)
\right) = 0, \ \ (x,y) \in \Omega \ .
\end{eqnarray}
If one applies a voltage $f$ on the boundary of $\Omega$, one then finds a unique solution to (\ref{eq:Conduct1}) together with the Dirichlet boundary condition
\begin{eqnarray}
u(x,y) = f(x,y), \ \ (x,y) \in \partial \Omega \ .
\end{eqnarray}
On the other hand, if one knows the current density on the boundary (a function $g(x,y)$), then one can also find a solution to the equation (\ref{eq:Conduct1}) together with the {\it Neumann} boundary condition
\begin{eqnarray}
\gamma(x,y) \frac{\partial u}{\partial \nu} (x,y) = g(x,y), \ (x,y) \in \partial \Omega \ , 
\end{eqnarray}
where $\nu$ is the outward normal to $\partial \Omega$.  The business at hand centers on the Dirichlet-to-Neumann mapping, which produces an appropriate Neumann data function $g$ given Dirichlet data $f$.

The inverse problem, in a nutshell, is to start with the full knowledge of the Dirichlet-to-Neumann mapping, and to determine the conductivity $\gamma(x,y)$ everywhere in $\Omega$.  The steps, as outlined in \cite{IMNS2004}, are as follows:
\begin{enumerate}
\item[A.] From the Dirichlet-to-Neumann mapping, construct the so-called scattering transform $t(k)$.  The steps in constructing this {\it without direct use of the conductivity equation} are outlined, for example, in \cite{HamHerMueVHer2012}.
\item[B.] Re-cast the original conductivity equation as a Schr\"{o}dinger equation in the plane by writing $q = \gamma^{-1/2} \Delta \gamma^{1/2}$, and $\tilde{u} = \gamma^{1/2} u$, so that
\begin{eqnarray}
\left( - \Delta + q \right) \tilde{u} = 0 \ ,\  (x,y) \in \Omega \ ,
\end{eqnarray}
and then extend $q$ to the entire plane, and seek a solution $\tilde{u} = \psi(x,y,k)$ with the following asymptotics as $|z| \to \infty$:
\begin{eqnarray}
e^{- i k z } \psi(x,y,k) = 1 + \mathcal{O} \left( \frac{1}{z} 
\right) \ .
\end{eqnarray}
\item[C.] The quantity $\psi(x,y,k)$ solves a D-bar equation {\it in the variable } $k$, which is best expressed in terms of $\mu = e^{- i k z } \psi(x,y,k) $:
\begin{eqnarray}
\dbar \mu = \frac{t(k)}{4 \pi \overline{k} } e^{2 i (k_{1} x - k_{2} y )} \overline{\mu} \ .
\end{eqnarray}
\item[D.]  The conductivity is extracted from $\mu$ by evaluating it at $k=0$:
\begin{eqnarray}
\lim_{k \to 0} \mu(x,y,k) = \gamma(x,y)^{1/2} \ .
\end{eqnarray}

\end{enumerate}

\subsection{Application 2:  Integrable nonlinear partial differential equations in $2+1$ dimensions} 
\label{Sec:DSII} The defocusing Davey-Stewartson II equation is the following nonlinear partial differential equation
\begin{eqnarray}
\label{eq:DSIIa}
&&i q_{t} + \frac{1}{2} \left( q_{xx} - q_{yy} \right) = - |q|^{2} q + \varphi q \\
\label{eq:DSIIb}
&& \varphi_{xx} + \varphi_{yy} = -2 \left( |q|^{2} \right)_{xx}.
\end{eqnarray}
In the above, $q = q(x,y,t) = q(z,t)$ is complex-valued, and subscripts denote partial derivatives.  Solutions of this equation can be interpreted (in an idealized sense) as representing the envelope of a two-dimensional surface wave propagating unidirectionally with a specified wave number.

{\bf Comment on notation}   Here and in what follows we will  write $\psi(z,k)$ even when $\psi$ is not analytic in $z$ or $k$.  

The Davey-Stewartson II equation (\ref{eq:DSIIa} - \ref{eq:DSIIb}) is completely integrable, possessing a Lax-pair
\begin{eqnarray}
\label{eq:specprob}
&&\psi_{x} + i \sigma_{3} \psi_{y} = \pmtwo{0}{q}{\overline{q}}{0} \psi, \\
\label{eq:TimeEvolve}
&&\psi_{t} = \pmtwo
{i \dbar^{-1} \partial \left( |q|^{2} \right)/2}{-i \partial q}
{i \dbar \overline{q}} {-i \dbar^{-1} \partial \left( |q|^{2} \right)/2} \psi - \pmtwo{0}{q}{\overline{q}}{0} \psi_{y} + i \sigma_{3}\psi_{yy} \ ,
\end{eqnarray}

The equation (\ref{eq:specprob}) is referred to, by analogy with previous 1-dimensional cases, as the spectral problem associated to the DS-II equation.  It can equivalently be written as
\begin{eqnarray}
\label{eq:specprob2}
\pmtwo{\dbar}{0}{0}{\partial} \psi = \frac{1}{2}\pmtwo{0}{q}{\overline{q}}{0} \psi \ .
\end{eqnarray}
The {\it direct problem} is then to seek a vector-valued solution $\psi = \psi(z,k) = \left( \begin{array}{c}
\psi_{1} \\
\psi_{2} \\ \end{array}
\right)$ with the following asymptotic behavior as $|z| \to \infty$:
\begin{eqnarray}
\label{eq:specprob2norm1}
&&\lim_{|z| \to \infty} \psi_{1} e^{-kz} = 1, \\
\label{eq:specprob2norm2}
&& \lim_{|z| \to \infty} \psi_{2} e^{-\overline{k} \overline{z}} = 0 \ .
\end{eqnarray}
As mentioned above, the quantity $\psi$ is referred to as a CGO solution.  The {\it reflection coefficient}, $r(k)$, is encoded in the sub-leading term in the asymptotic expansion of $\psi$, via
\begin{eqnarray}
\psi_{2} e^{ - \overline{k} \overline{z}} = \frac{r(k)}{\overline{z}} + \mathcal{O} \left( \frac{1}{|z|^{2}} \right) \ .
\end{eqnarray}
It is customary (though not necessary) to re-cast the equations (\ref{eq:specprob2}) as a D-bar problem by defining
\begin{eqnarray}
\mu_{1} = \psi_{1} e^{ - k z}, \ \ \ \ \mu_{2} = \overline{\psi}_{2} e^{- k z} \ ,
\end{eqnarray}
for then $\mu = \left( \begin{array}{c} \mu_{1} \\ \mu_{2} \\ \end{array} \right)$ solves the D-bar problem
\begin{eqnarray}
\dbar \mu = \frac{qe^{ \overline{k}\overline{z} - kz}}{2} \pmtwo{0}{1}{1}{0} \overline{\mu} \ .
\end{eqnarray}
Diagonalizing the above system produces equations of the form (\ref{eq:001}).  Indeed, the matrix $\pmtwo{0}{1}{1}{0}$ has eigenvalues $\pm 1$, and can be diagonalized via
\begin{eqnarray}
\pmtwo{0}{1}{1}{0} = \pmtwo{1}{1}{1}{-1} \sigma_{3} \left( 
\frac{1}{2} \pmtwo{1}{1}{1}{-1} \right) \ .
\end{eqnarray}
If one sets
\begin{eqnarray}
M = \frac{1}{2} \pmtwo{1}{1}{1}{-1} \mu, 
\end{eqnarray}
then the entries of $M = \left( \begin{array}{c} M_{1} \\ M_{2} \\ \end{array} \right)$ satisfy
\begin{eqnarray}
\label{eq:021}
&&\dbar M_{1}  = \frac{qe^{ \overline{k}\overline{z} - kz}}{2} \overline{M_{1}} , \\
\label{eq:022}
&&\dbar M_{2}  = -\frac{qe^{ \overline{k}\overline{z} - kz}}{2} \overline{M_{2}}, 
\end{eqnarray}
two equations of the form (\ref{eq:001}), both quantities normalized so that
\begin{eqnarray}
\label{eq:022b}
M_{j} = 1 + \mathcal{O}\left( \frac{1}{z}\right) \ \mbox{ as } |z| \to \infty \ .
\end{eqnarray}.  
The reflection coefficient is obtained from (\ref{eq:021})-(\ref{eq:022}) via
\begin{eqnarray}
\label{eq:023}
r=\frac{1}{2} \left( M_{1}^{(1)} - M_{2}^{(1)} \right)
\end{eqnarray}
where $M_{1}^{(1)}$ and $M_{2}^{(1)}$ represent the $\frac{1}{z}$ coefficients of $M_{1}$ and $M_{2}$ as $z \to \infty$:
\begin{eqnarray*}
M_{1} = 1 + \frac{M_{1}^{(1)}}{z} + \cdots, \ \ \ \ \ M_{2} = 1 + \frac{M_{2}^{(1)}}{z} + \cdots \ .
\end{eqnarray*}

As $q = q(x,y,t)$ evolves according to (\ref{eq:DSIIa}-\ref{eq:DSIIb}), the reflection coefficient evolves as well, but the evolution of $r$ is quite simple:
\begin{eqnarray}\label{rt}
r(k,t) = r(k,0) e^{ \frac{-it}{4}  \left( k^{2} + \overline{k}^{2} \right) } \ .
\end{eqnarray}
Better still, the reconstruction of the potential $q(x,y,t)$ given $r(k,0)$ is achieved via a D-bar problem {\it in the variable $k$}, so that in (\ref{eq:dbarrecA}-\ref{eq:dbarrecB}) below, $\dbar = \frac{\partial}{\partial \overline{k} }$.  Already placing things in diagonal form, the problem is:
\begin{eqnarray}
\label{eq:dbarrecA}
&&\dbar \tilde{M}_{1}  = e^{ \overline{k}\overline{z} - kz} \overline{r(k)} \ \overline{\tilde{M}_{1}} \ , \\
\label{eq:dbarrecB}
&&\dbar \tilde{M}_{2} = -e^{ \overline{k}\overline{z} - kz} \overline{r(k)}\  \overline{\tilde{M}_{2}}\ , 
\end{eqnarray}
together with the normalization conditions 
\begin{eqnarray}
&& \tilde{M}_{1} = 1 + \frac{\tilde{M}_{1}^{(1)}}{k} + \mathcal{O} \left( \frac{1}{|k|^{2}} \right) \ \mbox{ as } |k| \to \infty, \\
&& M_{2} = 1 + \frac{\tilde{M}_{2}^{(1)}}{k} + \mathcal{O} \left( \frac{1}{|k|^{2}} \right) \ \mbox{ as } |k| \to \infty \ .
\end{eqnarray}
The potential $q(x,y,t)$ is then obtained via
\begin{eqnarray}
\label{eq:qRec}
q(x,y,t) = \tilde{M}_{1}^{(1)} - \tilde{M}_{2}^{(1)} \ .
\end{eqnarray}

\subsection{Application 3:  2D orthogonal polynomials and Normal Matrix Models in Random Matrix Theory}
The orthogonal polynomials in this example are denoted $\{ p_{j} = p_{j}(z) \}_{j=0}^{\infty}$, $z$ is the usual complex variable, and the orthogonality is with respect to the two-dimensional measure $W(x,y) dA$:
\begin{eqnarray}
\iint_{\mathbb{R}^{2}} \overline{p}_{k} p_{k} W(x,y) dA = \delta_{jk} \ , \ \ j, k = 0, \ldots \ .
\end{eqnarray}
The orthogonality condition is equivalent to the conditions contained in the statement that the following solid Cauchy transform decays rapidly for $z \to \infty$:
\begin{eqnarray}
\iint \frac{ \overline{p}_{n}(z')}{z' - z} W(x',y') dA' = \mathcal{O} \left(
z^{-(n+1)} \right)
\ \mbox{ as } z \to \infty \ .
\end{eqnarray}
Putting this together, one sees that the row vector 
\begin{eqnarray}
\vec{\mu} :=\left( p_{n} , \ \ \frac{1}{\pi} \iint \frac{ \overline{p}_{n}(z')}{z' - z} W(x',y') dA' \right)
\end{eqnarray}
solves a d-bar problem:  the row vector behaves as follows for $z \to \infty$:
\begin{eqnarray}
\label{eq:normInf}
\vec{\mu} \left(
\begin{array}{cc}
z^{-n} & 0 \\ 0 & z^{n} \\ \end{array}
\right) = \left( 1, 0 \right) + \mathcal{O}\left(\frac{1}{z^{n}}\right) \ ,
\end{eqnarray}
and $\vec{\mu}$ satisfies
\begin{eqnarray}
\overline{\partial} \vec{\mu} = \overline{\vec{\mu}} \left(
\begin{array}{cc}
0& -W(x,y) \\
0& 0 \\
\end{array}
\right) \ .
\end{eqnarray}

For applications to the Normal Matrix Model in Random Matrix Theory, the sort of weight functions that are of interest include those of the form $W = e^{- N\left(  |z|^{2} + u(x,y)\right)}$ where the function $u(x,y)$ is a harmonic function whose growth at $\infty$ is slower than quadratic, so that the integrals are finite.  While the analysis of orthogonal polynomials with respect to a weight supported on the real axis (or even, in some cases, supported on a contour) has been carried out via Riemann-Hilbert techniques, the asymptotic analysis of these polynomials (as the degree grows to $\infty$ with the parameter $N$) has defied analysis, except in a select few special examples.  On the other hand, their asymptotic behavior is central to a great many applications, from the theory of Normal Random Matrices through to approximation theory.

\subsection{Desiderata}  For the first two applications described above, what is really required from the D-bar problem is the solution {\it at a specific point}.  In the first example, the conductivity is obtained from the solution to the D-bar problem at $k=0$, while in the second example, the reflection coefficient is obtained from the coefficient of $\frac{1}{\overline{z}}$ in the expansion of the solution to the D-bar problem as $|z| \to \infty$.  (The third example, 2D orthogonal polynomials, is more demanding, as the solution is actually required (in applications to the Normal Matrix Model) for all $z \in \mathbb{C}$.

In all cases, however, reformulating the D-bar equation as an integral equation shows that a starting point is the existence theory for a solution $\psi$ in $L^{\infty}(\mathbb{C})$, which can then be promoted to a solution with greater regularity if desired.  

Thus, in the literature one finds existence theory in $L^{\infty}$, 
for potentials $Q$ in relatively weak spaces, and one of the driving quests has been to obtain existence and uniqueness results under essentially minimal assumptions on the potential $Q$.

\subsection{Brief description of the numerical method developed by Knudsen, Mueller, and Siltanen}

A numerical solution method for (\ref{eq:001}) was developed and implemented in \cite{KnMuSi2004}.  The approach taken was to express the fundamental equation (\ref{eq:001}), normalized so that $\psi = 1 + v$, with $v = \mathcal{O} \left( 1/|z|\right)$ as $|z| \to \infty$, in terms of the operator $\dbar^{-1}$, which is a well-known singular integral operator in $\mathbb{C}$.  The integral equation takes the form
\begin{eqnarray}
v = 1 + \frac{1}{\pi} \iint_{\mathbb{C}} \frac{Q(z')\overline{v}}{z-z'} \ d^{2}z' \ .
\end{eqnarray}
(Compare to equation (1.3) in \cite{KnMuSi2004}, with $k$ replaced by $z$ and $T(k')$ replaced by $Q(z')$.)  Because of the application to the inverse conductivity problem, they considered functions $Q$ of compact support.  This compact support lent itself to a periodic extension, which they showed restricts to the {\it actual solution} on the original support.  Casting the problem in a periodic setting in turn led them to the use of the fast-Fourier transform in order to evaluate the discretization of the integral operator
\begin{eqnarray}
\mathcal{I}(\varphi)=
\int_{-s}^{s}
\int_{-s}^{s} g(z - z') Q(z')  \varphi(z') dz_{1} dz_{2} \ ,
\end{eqnarray}
where $g=g(z)$ represents the periodic extension of the function $\frac{1}{z}$ and the fundamental period domain is $[-s,s]^{2}$.

In discretizing and periodically extending the function $1/z$, the authors handled the singularity at $z=0$  by the relatively simple rule
\begin{eqnarray}
g_{h}(z) =  \left\{ \begin{array}{cc}
g(jh), & \mbox{ for } j \in \mathbb{Z}_{m}^{2} \setminus 0, \\
0, & \mbox{ for } j = 0 \ .\\
\end{array} \right.
\end{eqnarray}
See \cite{KnMuSi2004} for a description of the grid, and the finite integer lattice $\mathbb{Z}_{m}^{2}$.  

The authors showed that their method converges as the 
grid spacing tends to $0$, demonstrating that the method is a first 
order method.  The reason the method is first order and does not show 
spectral convergence of  a Fourier approach is this simple regularization: 
though a Riemann integral does not change if the integrand is just 
changed in one point, the same is not true for a spectral method, in 
particular not for a discrete Fourier transform. 
Changing the function in one point implies that it is not continuous, 
and it is well known that  the  Fourier coefficients $c_{n}$, $n\in 
\mathbb{Z}$ for a 
non-continuous function decrease as $1/n$. 
Starting from the successful regularization, 
the authors subsequently developed a two-grid method, based on the 
work of Vainikko \cite{Vainikko}. Note, however, that the authors 
consider functions of lower regularity than we, and for such 
functions spectral convergence cannot be observed even with the 
regularization to be discussed in the following sections. (Moreover, for the intended application to Electrial Impedance Tomography, the available data from measurements is sparse for this ill-posed problem, and so what is required is a lower order, fast approximation scheme.)

Prior to this work, there were other approaches to the numerical solution of the D-bar equation (\ref{eq:001}).  These approaches are described briefly in \cite{KnMuSi2004}, and we quote from that article, replacing citations with our own citation numbering scheme as needed (equation (1.1) in the quotation below is $\dbar v(k) = - T(k) \overline{v(k)} \ $ ):
\begin{quotation}
``The numerical solution of equation (1.1) was considered by Siltanen et al \cite{IsMuSi2000} for EIT in the numerical algorithm based on Nachman's uniqueness proof for the 2-D inverse conductivity problem \cite{Nach1996}. In \cite{IsMuSi2000, IsMuSi2001,MuSi2003,IsMuSi2002},
equation (1.1) was solved numerically by a 2-D adaptation of the method of product integrals \cite{Atkinson1989}. The numerical solution of equation (1.1) has also been
applied to EIT by Knudsen \cite{Knudsen2002} in the numerical algorithm based on the Brown-Uhlmann uniqueness proof for the 2-D inverse conductivity problem \cite{BrownUhl1997}. A fast, direct algorithm for the Lippmann-Schwinger equation in two
dimensions is found in \cite{Chen2002}.''
\end{quotation}
More details can be found in Chapters 14 and 15 of the text book \cite{MuSi2012}.  The basic approach to handling the inverse of the D-bar operator developed in \cite{KnMuSi2004} has been used in a number of different applications, see for example \cite{IMNS2004,KnLaMuSi2009,AsMuPaSi2010,KoLaOlSi2013} and references cited therein.  Very recently, in \cite{IMNS2014}, a finite difference solver for the D-bar equation (which was shown to be second order) was implemented to reconstruct EIT images.

\section{Brief discussion of existence theory for equation (\ref{eq:001})}
\label{Sec:2}
In this section we present a brief summary of known existence results 
for solutions to D-bar problems. As mentioned above, the D-bar equation (\ref{eq:001}), re-written for the unknown $m$ (\ref{eq:038}) can be re-cast as a singular integral equation
\begin{eqnarray}
\label{eq:045}
m(x,y) = 1 + \frac{1}{ 2\pi} \iint_{\mathbb{C}} \frac{ Q(z') e^{ \overline{kz' } - kz'} \overline{m(z')}}{z - z'}d^{2} z' \ + \ \frac{1}{ 2\pi} \iint_{\mathbb{C}} \frac{ Q(z') e^{ \overline{kz' } - kz'} }{z - z'} d^{2}z' \ ,
\end{eqnarray}
an integral equation of the form
\begin{eqnarray}
\left( 1 - \mathcal{I}_{Q,k}\right) m = \mathcal{I}_{Q,k}(1) \ ,
\end{eqnarray}
where
\begin{eqnarray}
\mathcal{I}_{Q,k}(f) = \frac{1}{2 \pi} \iint_{\mathbb{C}} 
\frac{ Q(z') e^{ \overline{kz' } - kz'}  \ \overline{f(z')} \ }{z - z'} d^{2} z' \ .
\end{eqnarray}
The existence theory for this equation has been explained in a number of places;  we refer to \cite{Sung1992} and \cite{Perry2012} and references contained therein for the following description of the existence theory.  (In those references, the authors were dedicated to establishing global existence for solutions of the Davey-Stewartson II initial value problem.  Along the way they had to tackle the existence, uniqueness, and regularity of a coupled version of the above D-bar problem.)

The operator $\mathcal{I}_{Q,k}$ is first shown to be a Fredholm operator of index zero.  The goal is to prove that there is nothing in the kernel of the operator $1 - \mathcal{I}_{Q,k}$.  If this is successful, Fredholm theory implies that $1 - \mathcal{I}_{Q,k}$ is an invertible operator.  A uniform bound on the {\it inverse} is obtained by first establishing a bound for large $k$, and then using continuity in $k$ for the bounded subset of the $k$ plane remaining.  

This is the original approach taken by Sung \cite{Sung1992}, who established the existence in $L^{\infty}(\mathbb{R}^{2})$, with potentials in $L^{1,\infty}(\mathbb{R}^{2})$ (the space of complex valued bounded and Lebesgue integrable functions on $\mathbb{R}^{2}$), and then proceeded to establish differentiability properties of the solution when the potential has derivatives in $L^{1,\infty}(\mathbb{R}^{2})$.  One particular result of his work (but not the only one) is that the {\it reflection coefficient} $r(k_{1}, k_{2})$ is in the Schwartz class $\mathbb{S}_{k_{1},k_{2}}\left( \mathbb{R}^{2} \right) $ of functions if the potential $Q(x,y)$ is in $\mathbb{S}_{x,y}\left( \mathbb{R}^{2} \right)$.

Perry \cite{Perry2012} iterated the integral equation one time, and thus analyzed integral operators of the form
\begin{eqnarray}\label{perryform}
\mathcal{T}(\psi) = \iint_{\mathbb{C}}\frac{e^{ \overline{kz'}-kz'}Q}{z-z'}\iint_{\mathbb{C}} \frac{e^{kz''-\overline{kz''}}\overline{Q}}{z'-z''}\psi(z'') d^{2}z'' \ d^{2} z' \ .
\end{eqnarray}
He then follows the basic Fredholm theory approach outlined above, 
assuming now that the potential $Q$ is in $H^{1,1}(\mathbb{R}^{2}) = 
\left\{f\in L^{2}(\mathbb{R}^{2}): \ \nabla f, \left( 1 + |\cdot| 
\right)f(\cdot) \in L^{2}\right\}$.  Note that we also use a once 
iterated version of the integral equation (\ref{eq:MainEQa}) in our 
numerical approach, but there the iteration is done in Fourier space. 

Perry proves that if $Q \in H^{1,1}(\mathbb{R}^{2})$, then the 
integral equation he considers (a coupled version of (\ref{eq:045})) has a unique solution in $L^{\infty}(\mathbb{R}^{2})$ which in addition is H\"{o}lder continuous in $z$ of order $\alpha$ for any $\alpha \in (0,1)$.  

Since the main effort is in establishing that the operator $( 1 - 
\mathcal{I}_{Q,k})$ is invertible, this implies immediately that {\it our equation} (\ref{eq:045}) possesses a solution in $L^{\infty}(\mathbb{R}^{2})$ which is also H\"{o}lder continuous with exponent $\alpha$ for any $\alpha \in (0,1)$, and $m$ actually converges to $1$ as $|z| \to \infty$.

Moreover, if we assume that our potential is in $\mathbb{S}_{x,y} 
\left( \mathbb{R}^{2}\right)$, then once a solution exists which is both bounded and H\"{o}lder continuous for $z \in \mathbb{C}$, standard potential theory (see, for example, \cite[pp. 53-61]{GilTru}) implies that the right hand side of (\ref{eq:045}) is H\"{o}lder continuous with exponent $\alpha + 1$, which clearly implies that $m$ is infinitely differentiable, and possesses a complete expansion in powers of $\frac{1}{z}$ for $|z| \to \infty$.  

We summarize these facts in a Theorem, which is attributable to the above collection of references:
\begin{theorem}\label{thrm:01}
Suppose that $Q \in \mathbb{S}_{x,y}(\mathbb{R}^{2})$.  Then the solution $m$ to the $\dbar$ equation (\ref{eq:039}) exists, possesses an expansion in powers of $\frac{1}{z}$ for $z \to \infty$, which may be differentiated term-by-term to obtain an expansion for derivatives of $m$.  Finally,
\begin{eqnarray}
\label{eq:049}
\dbar m \in \mathbb{S}_{x,y}(\mathbb{R}^{2}) \ .
\end{eqnarray}
\end{theorem}
The last conclusion, (\ref{eq:049}), follows from the $\dbar$ equation itself, and is the basis for a spectrally based numerical method.  Indeed, it provides rigorous justification for equation (\ref{eq:08}), re-written here for convenience (recall that $S(\xi) = \xi \hat{m}(\xi) )$.

\begin{eqnarray*}
(\ref{eq:MainEQa})\hspace{2.0in}
S(\xi) =  - i \mathcal{F}\left\{  Q e^{ \overline{kz}-kz} \left(\overline{ \mathcal{F}^{-1} \left( \frac{1}{\xi} S(\xi) \right) }  \right) \right\} 
 -i  \mathcal{F}\left\{  Q e^{ \overline{kz}-kz}  \right\} \ . \ \ \ \ \ \ \ \ \ 
\end{eqnarray*}

We close this section with a useful fact.
\begin{lemma}  Suppose that $h \in \mathbb{S}_{x,y}(\mathbb{R}^{2})$.  Then 
\begin{eqnarray*}
\dbar^{-1}(h) := \mathcal{F}^{-1}  \left( \frac{-2i}{\xi} \mathcal{F}(h) \right)
\end{eqnarray*}
is an infinitely differentiable function, possessing a complete expansion in powers of $\frac{1}{z}$ as $|z| \to \infty$.  
\end{lemma}
One may see this by writing
\begin{eqnarray}
\mathcal{F}^{-1} \left( \frac{-2i}{\xi} h(\xi_{1},\xi_{2}) \right) = 
\frac{-i }{\pi } \int_{\mathbb{R}^{2}} h(r \cos{\theta}, r 
\sin{\theta} )  e^{ i r \left( x \cos{\theta} + y \sin{\theta} 
\right)-i\theta} 
d \theta d r \ ,
\end{eqnarray}
which is clearly infinitely differentiable in $x$ and $y$.  The expansion in powers of $\frac{1}{z}$ follows either from a stationary phase analysis of this integral, or from the alternative representation of $\dbar^{-1}(h)$:
\begin{eqnarray}
\dbar^{-1}(h) = \frac{1}{\pi} \iint_{\mathbb{R}^{2}} \frac{h(x',y')}{z - z'} d^{2} z' \ .
\end{eqnarray}

\section{Numerical algorithm for the solution of the D-bar problem}
\label{Sec:3}
In this section we present an algorithm for the numerical solution of 
the integral equation (\ref{eq:MainEQa}) for small $k$. We will take $k=0$, but the algorithm works for small $k$ as well.  The approach 
follows in principle \cite{KnMuSi2004}, but uses a more sophisticated 
regularization of the singular integrand in 
 (\ref{eq:MainEQa}) to provide a sufficiently 
smooth integrand for the numerical computation of the integral. 

\subsection{Fourier approach and GMRES}

It is a well known fact from Fourier analysis that the Fourier 
transform of a Schwartz function is itself in a Schwartz space. This is the 
analytic basis of the efficient implementation of so-called spectral methods, in 
our case discrete Fourier series, which are used in this paper as 
follows: the Schwartz functions to be studied are considered on a 
sufficiently large computational domain,
\begin{equation}
    x\in L_{x}[-\pi,\pi], \quad y\in L_{y}[-\pi,\pi],
    \label{domain}
\end{equation}
that the function and its 
derivatives decrease to machine precision (here roughly $10^{-16}$, but in 
practice limited due to rounding errors to $10^{-14}$ or less) at the boundaries. This 
allows to treat Schwartz functions as smooth periodic functions 
within machine precision. In this case the Fourier coefficients 
decrease also exponentially. Approximating the Fourier series with a 
discrete Fourier series (computed efficiently via an FFT),
which corresponds in a sense to a 
truncation, the numerical error by neglecting the Fourier 
coefficients $c_{n}$ with $n>N/2$ will decrease exponentially with 
$N$, i.e., faster than any power of $1/N$. A method being in this sense of infinite order is called a 
spectral method. We always choose the wave numbers in the FFT as 
\begin{equation}
    \xi_{1}\in[-N_{x}/2+1,N_{x}/2]/L_{x},\quad 
    \xi_{2}\in[-N_{y}/2+1,N_{y}/2]/L_{y},
    \label{xidomain}
\end{equation}
where $N_{x}$ and $N_{y}$ are the number of Fourier modes in the $x$ 
and $y$ direction respectively. The values of $N_{x}$ and $N_{y}$ are 
chosen in a way that the Fourier transform of the studied Schwartz 
function decreases to machine precision as $\xi_{1}$ and $\xi_{2}$ approach the boundary of the computational domain. 

The discrete Fourier approach implies that 
the integral equation (\ref{eq:MainEQa}) is approximated by a finite dimensional system of 
linear equations. Writing the $N_{x}\times N_{y}$ matrix $S$ as a vector with 
$N_{x}N_{y}$ components, this system of equations has the form 
$\mathcal{A}S=b$. As in \cite{KnMuSi2004}, we use GMRES 
\cite{gmres} Krylov subspace 
techniques to solve this linear system. The advantage of GMRES is 
that the $N_{x}N_{y}\times N_{x}N_{y}$ matrix $\mathcal{A}$ does not 
have to be stored (this would make the use of parallel computers 
necessary for some of the values of $N_{x}$, $N_{y}$ used here), just 
its action on a vector has to be known.   The price to pay for this 
memory efficient implementation is that the equation $\mathcal{A}S=b$ 
is solved iteratively with a certain tolerance that can be chosen. 
For the high precision experiments to be presented below, we 
typically use a tolerance of the order of machine precision (in 
practical application, this can be of course chosen larger). In 
practice GMRES converges rapidly in the studied examples.   

\subsection{Regularization of the integrand}
To obtain spectral convergence the goal is to provide smooth periodic functions for all numerical evaluations involving an 
FFT.  Obviously this is not the case for the 
function $S(\xi)/\xi$ in (\ref{eq:MainEQa}), which is not even bounded 
for $\xi\to0$, though $S(\xi)$ is a smooth function, as explained in Theorem \ref{thrm:01}. 
To address the problem of $S(\xi)/\xi$ not being in the Schwartz class, the idea is to write the 
related integral in the form 
\begin{equation}
 \mathcal{F}^{-1} \left( \frac{1}{\xi} S(\xi) \right) 
= \mathcal{F}^{-1} \left( \frac{1}{\xi} 
    (S(\xi)-G(\xi) \right) 
    + 
\mathcal{F}^{-1} \left( \frac{1}{\xi} G(\xi) \right) 
    \label{Sreg},
\end{equation}
where $G(\xi)$ is chosen in a way that $(S(\xi)-G(\xi))/\xi$ is 
(within numerical precision) a smooth rapidly decreasing function, 
and that $\mathcal{F}^{-1} \left( G(\xi)/\xi \right) $ can be given 
in explicit form analytically up to constants which can also be computed with spectral 
accuracy. If this is possible, spectral convergence can be expected 
since $\mathcal{F}^{-1} \left( G(\xi)/\xi \right) $ will be 
multiplied with a Schwartz function, and thus the Fourier transform 
on the right hand side of (\ref{Sreg}) can be computed with the 
wanted precision.

To identify a suitable function $G(\xi)$ in (\ref{Sreg}), we use the 
following 
\begin{lemma}
    \label{lemma}
    The inverse Fourier transform of $\exp(-|\xi|^{2})/\xi$ is given 
    by
 \begin{eqnarray}
    \label{fourreg}
\mathcal{F}^{-1} \left( \frac{1}{\xi} e^{- |\xi|^{2} } \right) = \frac{i}{x + i y} \left(1 -  
e^{ - (x^2 + y^2)/4} 
\right)\ .
\end{eqnarray}
In addition one has 
\begin{align}
     \mathcal{F}^{-1} \left(\frac{ \overline{\xi}^{n} }{\xi}e^{ - |\xi 
    |^{2} } \right)&= 
    i(2i)^{n}\frac{n!}{z^{n+1}}\left\{1-\exp\left(-\frac{|z|^{2}}{4}\right)\sum_{k=0}^{n}
    \frac{1}{k!}\left(\frac{z\bar{z}}{4}\right)^{k}
    \right\}
    \label{fourreg5},
\end{align}
which is equivalent to 
\begin{align}
     \mathcal{F}^{-1} \left(\frac{ \overline{\xi}^{n} }{\xi}e^{ - |\xi 
    |^{2} } \right)&=i(2i)^{n}n!\left(\frac{\bar{z}}{4}\right)^{n+1}\sum_{k=0}^{\infty}
    \frac{1}{(k+n+1)!}\left(\frac{z\bar{z}}{4}\right)^{k}
    \label{fourreg5a}.
\end{align}
\end{lemma}
Relation (\ref{fourreg}) can be established by direct calculation, 
whereas the remaining equations follow from it via
\begin{equation}
    \mathcal{F}^{-1} \left(\frac{ \overline{\xi}^{n} }{\xi}e^{ - |\xi 
    |^{2} } \right) = 
    \left(-2i\partial_{z})^{n}\right)\frac{i}{z}\left(1-\exp(-|z|^{2}/4)\right)
    \label{fourreg3},
\end{equation}
by direct calculation. 

\begin{remark}
    The term in (\ref{fourreg5a}) is similar to so-called 
    $\phi$-functions appearing in ETD schemes, see for instance 
    \cite{KassT} and references therein.
    This formula is numerically problematic close to the origin (say 
    for $|z|<1$), because 
    of cancellation errors. There it 
     is more convenient to use the Taylor expansion 
     (\ref{fourreg5a}). This allows one to compute the function to 
     machine precision in the whole computational domain. 
\end{remark}

\begin{remark}
    Lemma (\ref{lemma}) permits the construction of a function $G(\xi)$ in 
    (\ref{Sreg}) such that $(S(\xi)-G(\xi))/\xi$ is a Schwartz 
    function. A possible choice is $G(\xi)=S(\xi)\exp(-|\xi|^{2})$, 
    which implies, however, that the last integral on the right hand side of 
    (\ref{Sreg}) can only be computed as a convolution,
    \begin{eqnarray}
    && \mathcal{F}^{-1} \left\{Se^{ - |\xi|^{2} }/\xi  \right\}  = \frac{1}{2\pi} \int \mathcal{F}^{-1}(S) 
    \frac{i}{z-z'} \left(  1 - e^{ - \overline{(z-z')} ( z - z')/4} \right)
     \ d^{2} z' \ .
    \end{eqnarray}
    This choice of $G(\xi)$ satisfies all requirements, but the 
    computation of the convolution is not very efficient since the 
    standard Fourier approach would involve functions of lower regularity 
    than wanted (exactly the function $\exp(-|\xi|^{2})/\xi$) 
    introduced since it has an analytically known  Fourier transform). Thus to 
    obtain a spectral  method, the convolution has to be computed 
    directly, which is computationally expensive, since the integrand is 
    an object in four real dimensions. It will be shown in this and 
    the
    following section how a two-dimensional approach can be realized. Nonetheless  we use a code 
    implementing this approach to provide an independent test of other codes.
\end{remark}

The nature of the singularity of $S(\xi)/\xi$ at $\xi=0$ is demonstrated by a two-dimensional Taylor expansion:
\begin{eqnarray*}
\frac{1}{\xi} S(\xi) = \frac{S(0,0)}{\xi} +  \sum_{j=1}^{M}\frac{1}{j!} \frac{\partial^{j}}{\partial \overline{x}^{j}} S(0,0) \left(\frac{\overline{\xi}^{j}}{\xi} \right)  + \mbox{regular} \ .
\end{eqnarray*}
The term ``$\mbox{regular}$'' possesses $M$ continuous derivatives, and the first two terms clearly possess only fewer derivatives at $\xi=0$.  Closer inspection shows that the first two terms can be thought of as the purely anti-holomorphic part of $S$, divided by $\xi$.  Thus we choose $G$ so as to eliminate these terms:
\begin{equation}
    G(\xi)=e^{-|\xi|^{2}}\sum_{n=0}^{M}\frac{1}{n!}\frac{\partial^{n}S(0)}{\partial \bar{\xi}^{n}}\bar{\xi}^{n}
    \label{Gchoice},
\end{equation}
where $M$ is taken large enough that the Fourier 
coefficients of $(S(\xi)-G(\xi))/\xi$ decrease to machine precision.  (To see this, observe that $(S(\xi)-G(\xi))/\xi$ possesses $M$ continuous derivatives, and decays rapidly to zero for $|\xi|$ approaching the boundary of the computational domain, and hence its Fourier coefficients will decrease rapidly as well).  The derivatives of $S$ in (\ref{Gchoice}) are numerically computed via
\begin{equation}
    \frac{\partial^{n}S(0)}{\partial \bar{\xi}^{n}} = 
    \left(-\frac{i}{2}\right)^{n}\mathcal{F}z^{n}\mathcal{F}^{-1}S
    \label{fourreeg4},
\end{equation}
i.e., also with FFT techniques which can be done with spectral 
accuracy since $S$ is a Schwartz function. Note, however, that for 
large $n$, the multiplication with $z^{n}$ in (\ref{fourreeg4}) will 
delimit the achievable numerical accuracy. 

\begin{figure}[htb!]
   \includegraphics[width=0.49\textwidth]{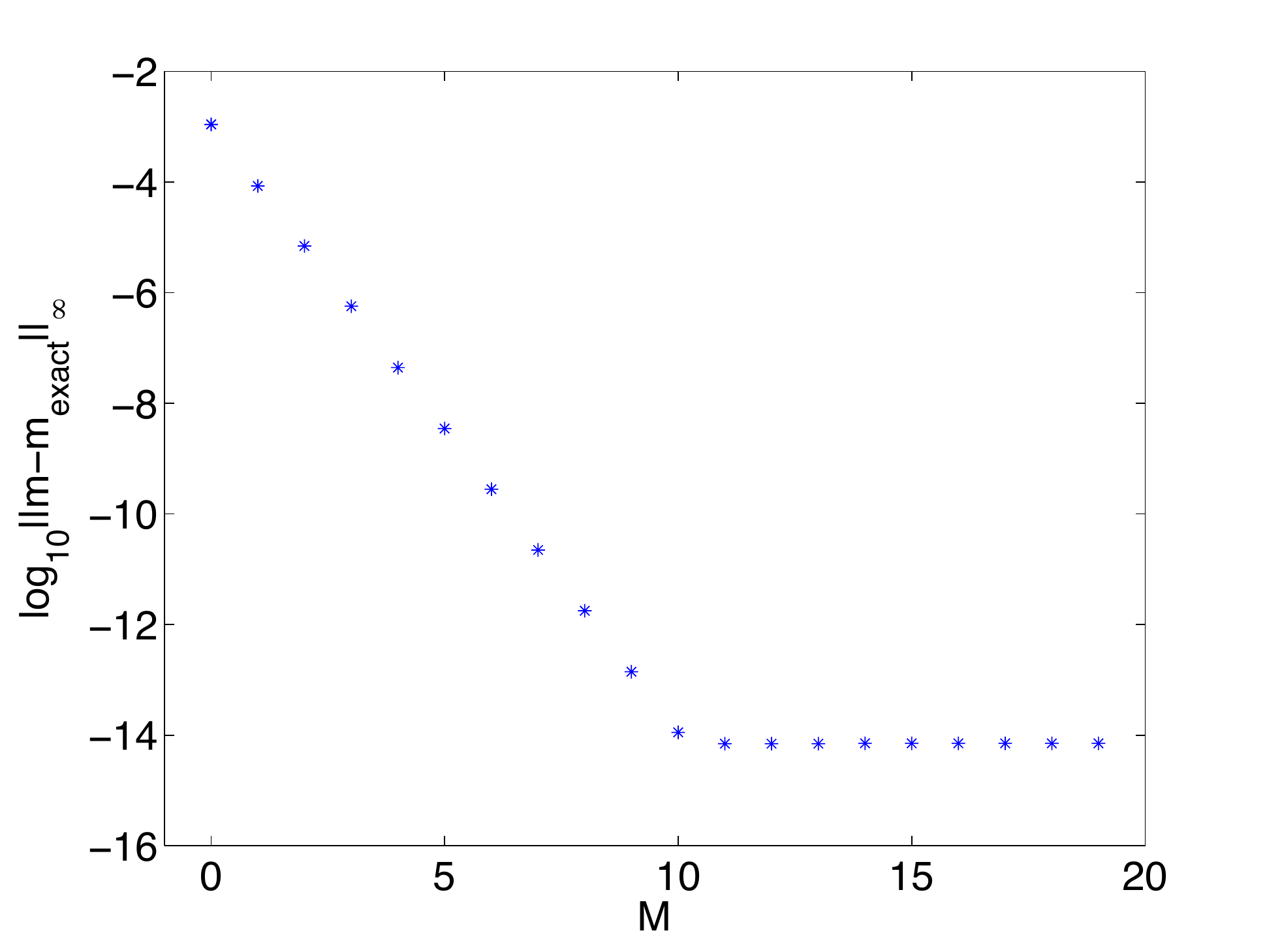}
   \includegraphics[width=0.49\textwidth]{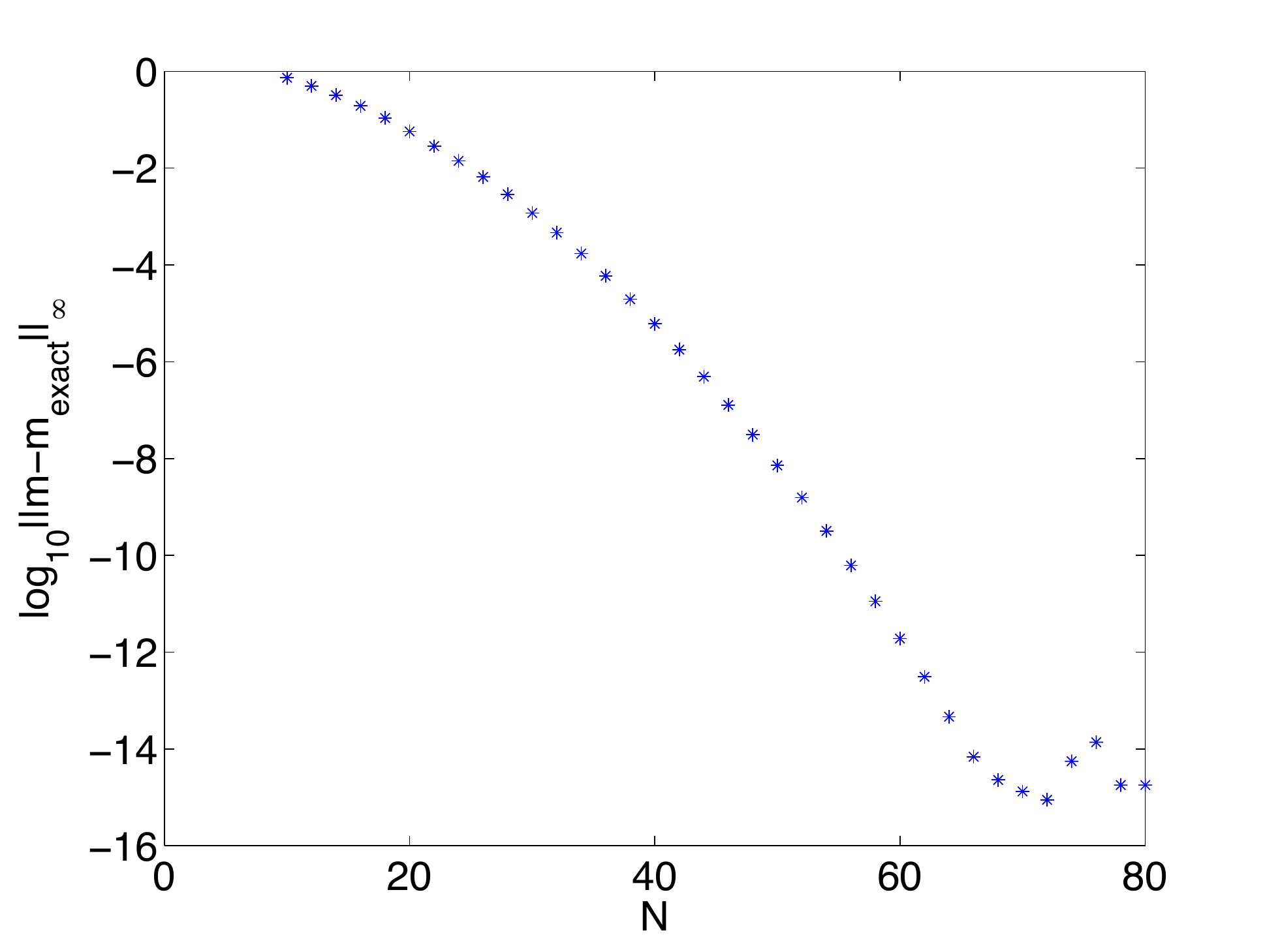}
 \caption{Difference between the numerical solution of the 
 D-bar inversion (\ref{dbartest}) with $a=1/2$, $b=1$ and $c=i$.  On 
 the left we plot the dependence on the number $M$ of terms in the Taylor expansion (\ref{Gchoice}), with $N=128$ Fourier modes in $x$ and $y$.  On the right we plot the dependence on the number of Fourier modes $N$ in $x$ and $y$, using $M=11$ terms in the Taylor expansion. }
 \label{fig:dbartest}
\end{figure}

To test the above approach to compute the inverse of the D-bar 
derivative for Schwartz functions via Fourier integrals, we consider 
the following test problem: It is straightforward to show by direct 
calculation that 
\begin{equation}
    \partial^{-1}_{\bar{z}}\exp(-a(z-b)(\bar{z}-c))=\frac{a}{z-b}\left(1-\exp(-a(z-b)(\bar{z}-c))\right)
    \label{dbartest},
\end{equation}
where $a>0$, $b$, $c$ are constants. In Fig.~\ref{fig:dbartest} we 
show how the difference between the numerical and the exact solution 
depends on the number of Fourier modes and the number $M$ of terms in the 
Taylor expansion of $S$ in (\ref{Gchoice}). For simplicity we use 
$N_{x}=N_{y}=N$ Fourier modes and $L_{x}=L_{y}=4$.

The left figure shows the dependence of the 
$L^{\infty}$ norm of the difference between numerical and exact 
solution on the number $M$ of terms in the Taylor 
series in (\ref{Gchoice}), using $N=128$ Fourier modes in both $x$ and $y$ directions.  It can be seen in the logarithmic plot 
that the numerical error decreases exponentially to machine 
precision. It saturates in this example roughly for $M=11$. This 
behavior can be understood as follows: as explained above, the function 
$(S(\xi)-G(\xi))/\xi$ possesses $M$ continuous derivatives and decreases rapidly to machine precision as $\xi$ approaches the boundary of the computational domain.  
It is well known that the 
Fourier coefficients of such a function function decrease algebraically 
in the number $N$ of Fourier modes as $N^{-(M+2)}$. This implies a 
numerical error of the same order if the series is approximated by a 
truncated series, and thus we should, and do, see an exponential decrease of the error with $M$. Note that due to the finite precision of the numerics, it does not make sense to 
take arbitrarily high values of $M$.  As mentioned above, numerical precision is 
lost for large values of $M$ in the computation of the derivatives 
(\ref{fourreeg4}) since the integrand is multiplied with high powers 
of $z$ in (\ref{fourreeg4}) and $\bar{\xi}$ in (\ref{Gchoice}). 
Still it is possible to use values of $M$ up to 20 since with larger 
$M$ also the functions $W_{n}$ become smaller 
at the boundaries, and   $\exp(-|\xi|^{2})$ is small there in any case. 

In a similar way we study in the right figure of 
Fig.~\ref{fig:dbartest} the dependence of the numerical error on the 
number of Fourier modes $N$ for fixed $M=11$ (the result does not 
change if a higher value of $M$ is taken). It can be seen that the 
numerical error also decreases essentially exponentially and reaches 
machine precision for $N\approx 2^{6}$.

\subsection{Algorithm for computing a solution to the D-bar equation 
(\ref{eq:MainEQa})
}
We present a brief summary for the algorithm to solve 
(\ref{eq:MainEQa}), with $k=0$.
\begin{enumerate}
    \item Introduce a discrete Fourier grid as in (\ref{domain}) and 
    (\ref{xidomain}).
\item Compute and store the functions 
\begin{eqnarray}\label{Wn}
\frac{1}{n!}W_{n} =   \frac{1}{n!}\mathcal{F}^{-1} \left(\frac{ \overline{\xi}^{n} }{\xi}e^{ - |\xi 
    |^{2} } \right) 
\end{eqnarray}
using definitions (\ref{fourreg5}) and (\ref{fourreg5a}).

\item Call GMRES with $b=-i\mathcal{F}Q$ and 
$\mathcal{A}S=S+i\mathcal{F}\left( Q \ \overline{\mathcal{F}^{-1}(S/\xi)} \right)$.
\item To compute $\mathcal{A}(S)$:
\begin{enumerate}
\item Compute the relevant Taylor coefficients of $S$ at $\xi = 0$, 
using FFT techniques:
\begin{eqnarray}\nonumber
n!c_{n}= \frac{\partial^{n}}{\partial \overline{\xi}^{n}}S (0) = \left( \frac{-i}{2} \right)^{n}\mathcal{F}\left(z^{n} \mathcal{F}^{-1}(S) \right)(0) \ .
\end{eqnarray}
\item Form
\begin{eqnarray}\nonumber
&&G(\xi) = e^{- | \xi |^{2} } \sum_{n=0}^{M} c_{n} \overline{\xi}^{n}, \\
&&\mathcal{F}^{-1} \left(  \frac{1}{\xi} G(\xi) \right)  = \sum_{n=0}^{M}c_{n} W_{n} \ .
\end{eqnarray}
\item Now compute the intermediate quantity:
\begin{eqnarray}\label{S0Taylor}
h=
\mathcal{F}^{-1} \left( \frac{S(\xi)}{\xi} \right) = \mathcal{F}^{-1} \left[
\frac{1}{\xi} \left( S(\xi) - G(\xi)  \right)
\right] \ + \ \sum_{n=0}^{M}c_{n} W_{n} \ .
\end{eqnarray}
\item And finally compute $\mathcal{A}(S)$:
\begin{eqnarray*}
\mathcal{A}(S) = S + i \mathcal{F} \left( Q   \ \overline{h} \  \right) \ .
\end{eqnarray*}
\end{enumerate}

\end{enumerate}
Note that no factorials $n!$ that could affect the accuracy 
need to be computed here,  since they cancel in the combination 
$W_{n}c_{n}$. Factorials just appear in the Taylor series 
(\ref{fourreg5a}). 
\begin{remark}
    The computational cost of the algorithm is mainly in the FFT 
    commands (always in two dimensions) in the GMRES iterations. The 
    functions $W_{n}$ in (\ref{Wn}) are computed and stored beforehand since they do 
    not change during the iteration. 
    The action of the matrix $\mathcal{A}$ in each GMRES iteration 
    requires two FFTs. In addition the computation of the Taylor 
    coefficients of $S$ in (\ref{S0Taylor}) requires one additional FFT (not two, since the 
    values are just needed for $\xi=0$, the second Fourier transform 
    in (\ref{fourreeg4}) can be replaced by a summation). To sum up, per 
    GMRES iteration 3 FFTs are necessary. Thus the computational cost 
    is higher than in \cite{KnMuSi2004} (one FFT more per iteration), 
    but the convergence of our approach is spectral instead of first 
    order. When dealing with potentials of sufficient regularity (as we do), we may reach much higher precision with less resolution.
\end{remark}

\section{Numerical computation of Complex Geometrical Optics solutions}
\label{Sec:4}

In this section we will present a numerical approach to efficiently 
solve the D-bar equation (\ref{eq:039}) for general $k\in\mathbb{C}$ 
with FFT methods. Since the 
equation (\ref{eq:MainEQa}) defined on $\mathbb{R}^{2}$ is 
approximated by a finite dimensional system on $\mathbb{T}^{2}$, the 
quality of the approximation will depend on the commensurability of the functions under consideration with the computational domain.  The 
presence of the term $e^{(\overline{kz}-kz)}$ in (\ref{eq:MainEQa}) 
implies a shift in Fourier space.  So, for example,  $\mathcal{F}(e^{(\overline{kz}-kz)}Q)$ will not 
vanish with numerical precision on the boundary of the computational 
Fourier domain for $|k|$ sufficiently large. Since this was the basis of the approach detailed in the previous section, a reformulation of (\ref{eq:MainEQa}) is 
necessary in order to apply the same techniques as before. 

\subsection{Dependence on the parameter $k$}
It is useful to re-express the fundamental equation (\ref{eq:MainEQa}) using an integral operator:
\begin{eqnarray}
\label{eq:MainEQDS}
S = \mathcal{K}(S) - i \mathcal{F}\left( Q e^{\overline{kz}-kz} \right) \ ,
\end{eqnarray}
where 
\begin{eqnarray}
\mathcal{K} \left(h\right) = -i \mathcal{F} \left( q e^{ \overline{kz} - kz} \left[\overline{ \mathcal{F}^{-1} \left( \frac{h}{\xi} \right)} \right] \right) \ .
\end{eqnarray}
The operator $\mathcal{K}$ is the composition of a simpler integral operator and a shift:
\begin{eqnarray}\nonumber
&&\mathcal{K} (h) = \mathcal{K}_{0} (h) \circ ( \xi + 2 i \overline{k}) \ , \\
&&\mathcal{K}_{0} = -i \mathcal{F} \left( q \left[\overline{ \mathcal{F}^{-1} \left( \frac{S}{\xi} \right) }\right] \right) \ .\label{shift}
\end{eqnarray}

For sufficiently small values of $k$,  the numerical approach 
described in the previous section provides a spectrally accurate discretization of the equation, and  the operator 
acting on $S_{\pm}$ can be directly inverted via 
GMRES \cite{gmres}. 

For larger values of $|k|$, this procedure encounters a new 
difficulty.  Asymptotic analysis of the integral operator 
$\mathcal{K}$ for $|k|$ large shows that it is a small norm operator, and we may invert the equation by Neumann series for $|k|$ sufficiently large.  The first term of this expansion is obviously the last term on the right hand side of (\ref{eq:MainEQDS}):
\begin{eqnarray}\label{S0}
S(\xi) \approx S^{(0)} = 
 - i  \mathcal{F}\left\{  Q e^{ \overline{kz}-kz}  \right\} = - i \hat{Q}(\xi + 2i \overline{k}) \ . 
\end{eqnarray}
Further iterations may be carried out, and so, for example, a more accurate approximation to the solution is obtained by substituting $S^{(0)}$ into the right hand side of equation (\ref{eq:MainEQDS}):
\begin{eqnarray}
S(\xi) \approx  - i \mathcal{F}\left\{  Q e^{\overline{kz}-kz} \left(\overline{ \mathcal{F}^{-1} \left( \frac{1}{\xi} \left( 
- i \hat{Q}( \xi + 2i \overline{k}) 
\right) \right) }  \right) \right\} 
 - i \hat{Q}( \xi +2 i \overline{k}) \ .
\end{eqnarray}
This may be simplified to the following form:
\begin{eqnarray}
S \approx S^{(0)} + S^{(1)},
\end{eqnarray}
where 
\begin{eqnarray}
\label{S1}
S^{(1)} = 
- i \mathcal{F} \left[
\frac{q}{ 2 (\overline{z} - \overline{k})} \left[
e^{ \overline{kz}-kz - |k|^{2}/4}- e^{ - |z|^{2}}
\right]
\right] \  \ .
\end{eqnarray}

\begin{figure}[htb!]
   \includegraphics[width=.49\textwidth]{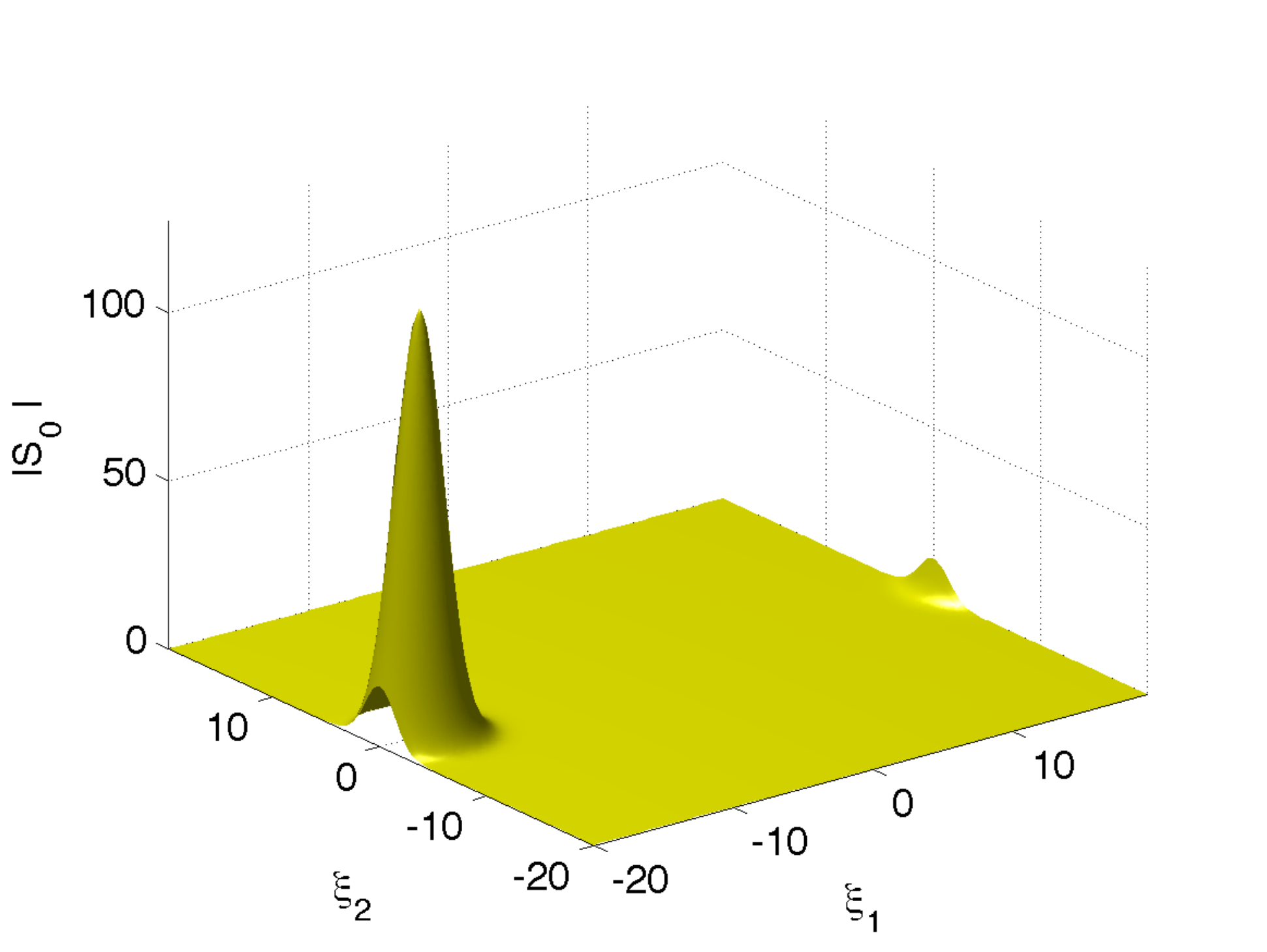}
   \includegraphics[width=.49\textwidth]{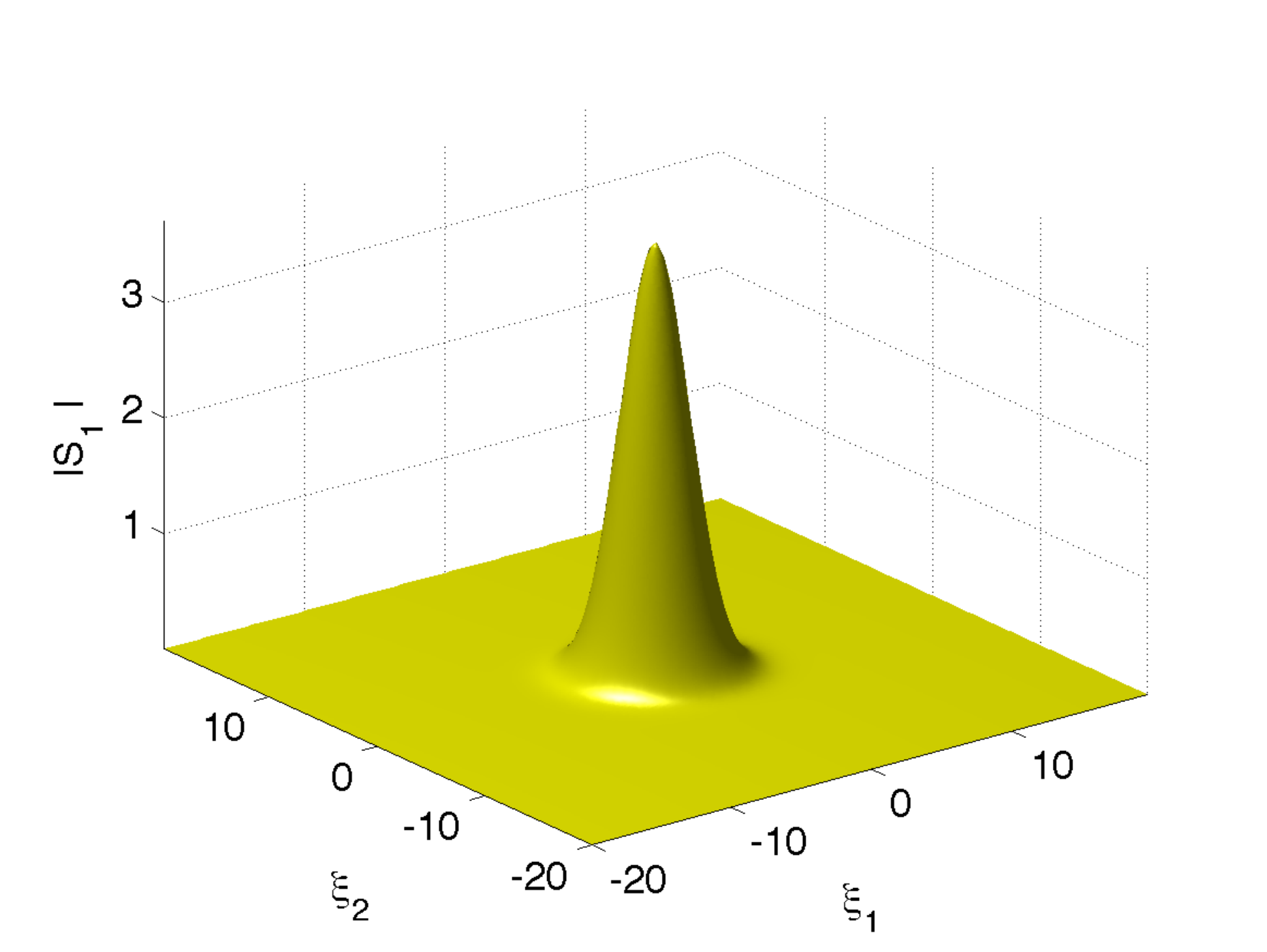}
 \caption{Modulus of the iterates $S_{0}$ (\ref{S0}) and $S_{1}$ (\ref{S1}) 
 for $k=8.5936$, $L=3.2$ and $N=128$. }
 \label{fig:kdepend}
\end{figure}

The above asymptotic description provides a key insight into the 
numerical issues we are facing for $k$ large.  The last term is a 
shift of $\hat{Q}$, and is supported in a neighborhood of $-2i 
\overline{k}$, which can be very near the boundary of the 
computational domain.  Moreover, the {\it first term} on the right 
hand side is itself the sum of two terms, the first also being 
supported in a neighborhood of $-2i \overline{k}$ (although 
exponentially suppressed by the term $e^{ - |k|^{2}/4}$), and the 
second term supported near $\xi = 0$. So the solution has support in two disparate regimes, one near the boundary of the computational domain and the other in the middle.
 This behavior can be seen 
for instance in Fig.~\ref{fig:kdepend}.

In the original numerical approach, we are dividing by $\xi$, which 
is not periodic, and not very small at the boundary of the 
computational domain.  This has the effect of introducing a 
discontinuity across the boundary of the computational domain, an 
anathema to spectral methods.  

However, if we know beforehand that a function $f$ has the property that its Fourier transform $\hat{f}$ has support centered at $-i \overline{k}$, we may compute $\dbar^{-1}(f)$ by exploiting the shift to our advantage:
\begin{eqnarray}
&&\dbar^{-1} (f) = -2i \mathcal{F}^{-1}\left[  \frac{\hat{f}}{\xi} \right] = 
-2i \mathcal{F}^{-1}\left[ \left( \frac{\hat{f}(\xi-2i\overline{k})}{\xi-2i\overline{k}} \right) 
\circ(\xi+2i \overline{k})
\right] \\
&& \ \ \ \ \ = -2 i e^{\overline{kz}-kz} \mathcal{F}^{-1}\left[ \left( \frac{\hat{f}(\xi-2i\overline{k})}{\xi-2i\overline{k}} \right) \right] \ .
\end{eqnarray}
Now in this last formula, the numerator, $\hat{f}(\xi - 2i \overline{k})$ is supported, by assumption, in a vicinity of $\xi=0$, and decays rapidly as $|\xi|$ approaches the boundary of the computational domain.  We may then apply the regularization approach of the previous section, but with the singularity at $\xi = 2i \overline{k}$ instead of $\xi = 0$.  Precisely, we take
\begin{eqnarray}
G(\xi) =  e^{-|\xi - 2i \overline{k}|^2} \ \sum_{n=0}^{M}\frac{1}{n!} \frac{\partial^{n}}{\partial \overline{\xi}^{n}}\hat{f}(0,0) (\overline{\xi}+2i k)^{n} \ , 
\end{eqnarray}
and then 
\begin{eqnarray}
\label{eq:dbarinvShift}
&&\dbar^{-1} (f) = -2 i e^{\overline{kz}-kz} \mathcal{F}^{-1}\left[ \left( \frac{\hat{f}(\xi-2i\overline{k}) - G(\xi) }{\xi-2i\overline{k}} \right) \right]  -2 i e^{\overline{kz}-kz} \mathcal{F}^{-1}\left[ \left( \frac{G(\xi)}{\xi-2i\overline{k}} \right) \right] \ .
\end{eqnarray} 

\begin{figure}[htb!]
   \includegraphics[width=0.49\textwidth]{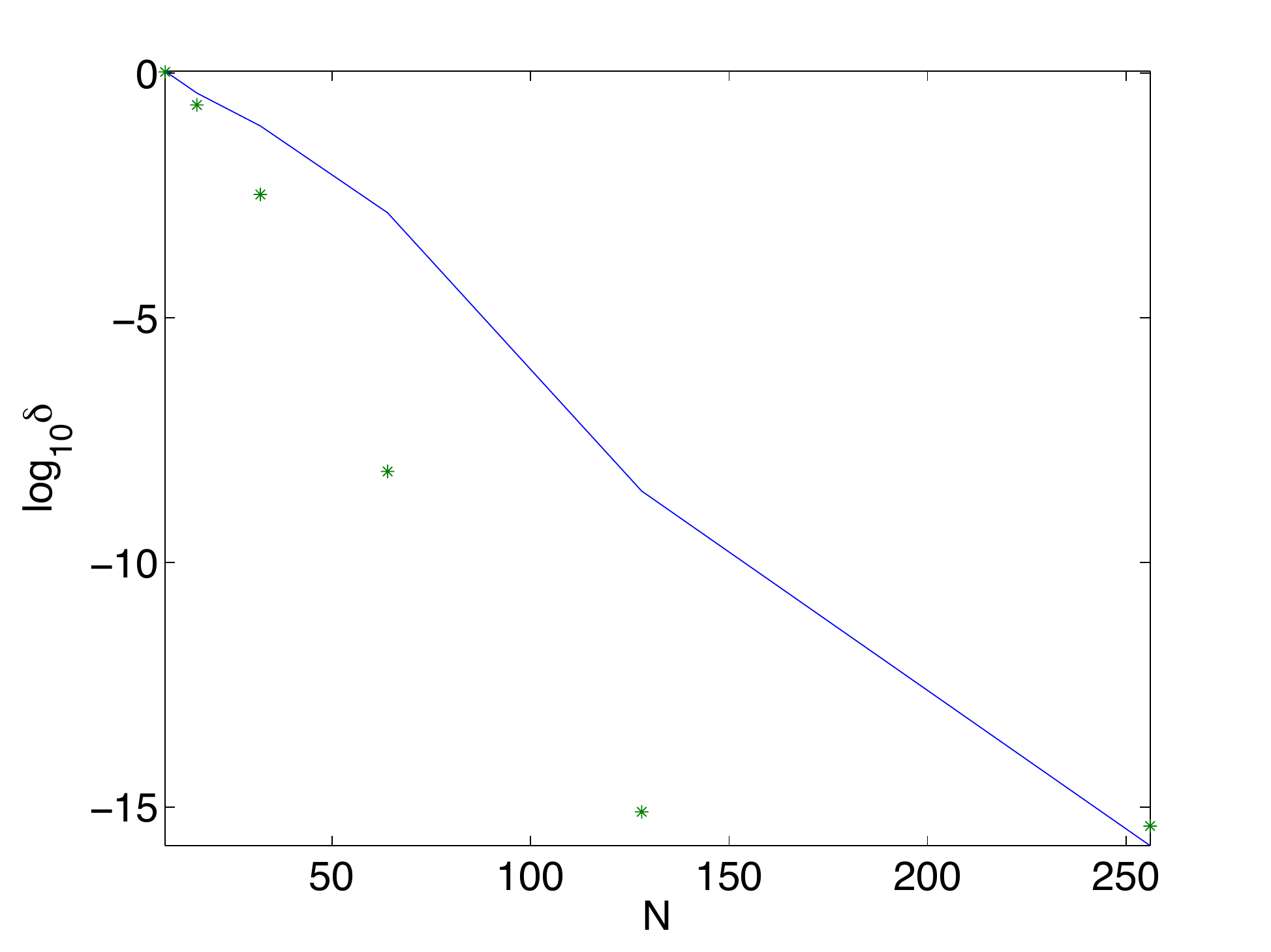}
   \includegraphics[width=0.49\textwidth]{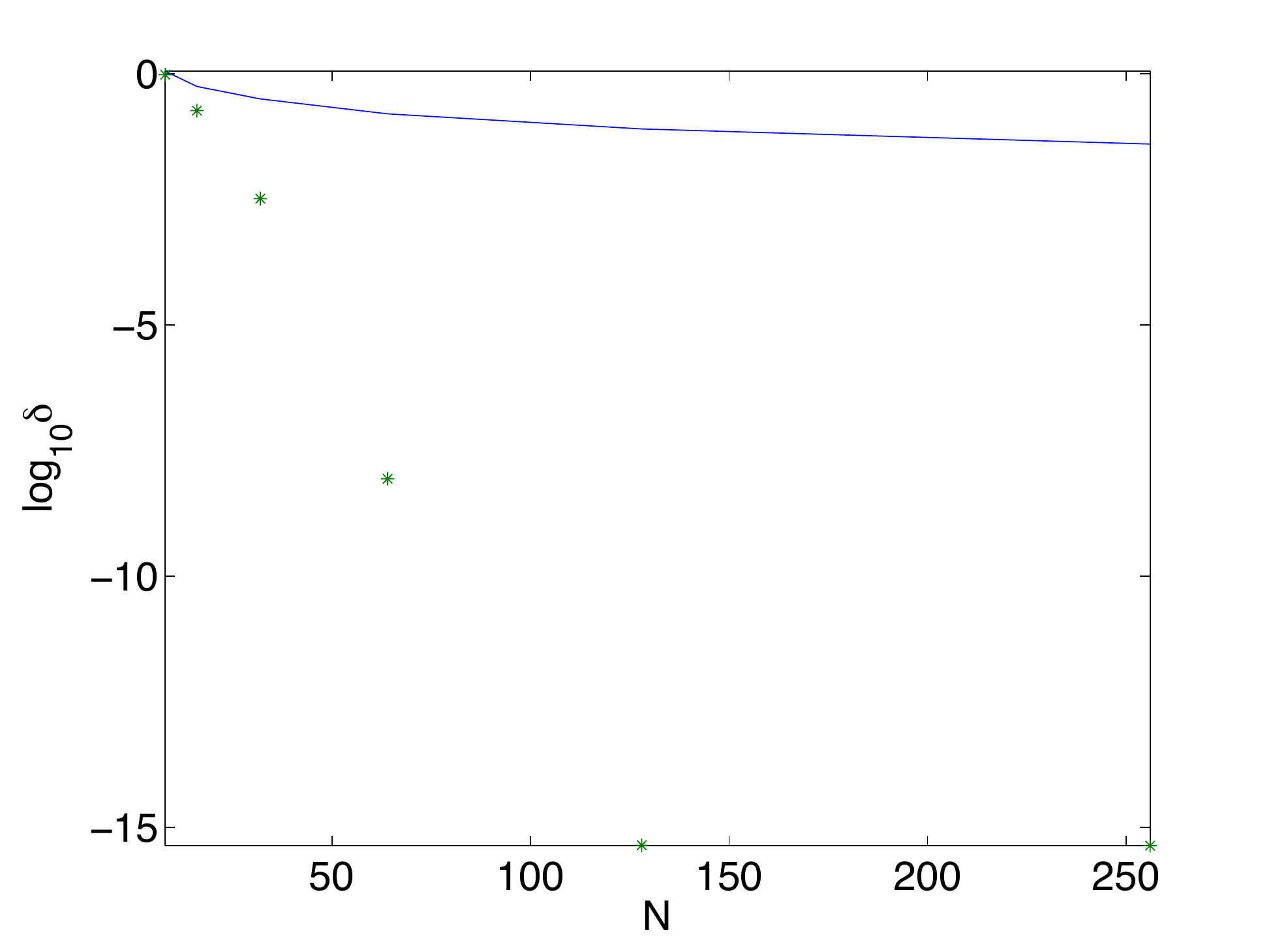}
 \caption{These are logarithmic plots of the dependence of the error in computing (\ref{dbartest}) on $N$, the number of Fourier modes.  The solid line shows the result of the computation without a shift in Fourier space, and the stars show the results with shifting.  On the left, the parameters in (\ref{dbartest}) were chosen so that the Fourier transform of (\ref{dbartest}) is supported mid-way between the origin and the boundary of the computational domain.  On the right the parameters were chosen so that the Fourier transform of  (\ref{dbartest}) is supported very near the boundary of the computational domain.  For these plots, the choice of $k$ varies with $N$, and the relevant parameter values are summarized in the Table shown in Fig.~\ref{Fig:FshiftTable}.}
 \label{fig:Fshift}
\end{figure}

We demonstrate the effectiveness of this shifting argument 
in Fig.~\ref{fig:Fshift} where we compute the solution to the example 
(\ref{dbartest}) for parameters leading to support of the function 
near the computational boundary with the approach of the previous 
section, both with 
and without the shifting argument (\ref{eq:dbarinvShift}). It can be seen that the 
numerical solution converges as expected very slowly to the exact 
one without a shift since the integrand is not small on the boundary. 
A shift of the integrand in Fourier space such that its maximum is 
localized near the center of the computational domain restores rapid 
convergence as in the previous section. 

However, it is important to recall that the CGO solutions which we seek are actually supported (for large $k$) in two disparate regions.
\begin{figure}[htb!]
\begin{tabular}{|c|c|c|c|}
\cline{1-4}
Fourier Modes & L & Left Figure $k$ & Right Figure $k$ \\
\hline
$2^{3}$ & $4.0$& $0.25(1+2i)$ & $0.5(1 + i)$  \\
\hline
$2^{4}$& $4.0$&$0.5(1+2i)$ & $1+i$\\
\hline
$2^{5}$ &$4.0$& $(1+2i)$&$2(1+i)$ \\
\hline
$2^{6}$& $4.0$&$2(1+2i)$&$4(1+i)$ \\
\hline
$2^{7}$& $4.0$&$4(1+2i)$ &$8(1+i)$ \\
\hline
$2^{8}$& $4.0$&$8(1+2i)$&$16(1+i)$\\
\hline
\end{tabular}
\caption{Parameter values for the plots shown in Fig.~\ref{fig:Fshift}.}  
\label{Fig:FshiftTable}
\end{figure}

\subsection{Iterated integral equation}

To capture the ``disparate regions of support" suggested by the asymptotic analysis, we  seek a solution $S$ of the form
\begin{eqnarray}
\label{eq:Sform}
S = f(\xi) + h\circ(\xi +2 i \overline{k}) \ ,
\end{eqnarray}
where $h$ and $f$ are functions supposed to vanish on the 
computational boundary and to have support in the vicinity of the 
origin. 
This leads to a pair of equations:
\begin{eqnarray}
&& h = \mathcal{K}_{0}(f) - i \mathcal{F}\left( q 
\right). \\
&&f = \mathcal{K}_{0}\left( h \circ(\xi +2 i 
\overline{k}) \right) \circ(\xi + 2i \overline{k})  \ ,
\end{eqnarray}

Replacing $f_{\pm}$ in the first equation above, we have
\begin{eqnarray}\label{eqh}
h- \mathcal{K}_{0}\left(
\mathcal{K}_{0}\left( h  \circ(\xi + 2i \overline{k}) \right) \circ(\xi + 2i \overline{k}) 
\right) =- i \mathcal{F}\left( q \right) \ .
\end{eqnarray}
Note that this equation is equivalent to Perry's integral equation 
(\ref{perryform}), but that we consider it in Fourier space. Then we find
\begin{eqnarray}
\label{eq:fpmdef}
f = \mathcal{K}_{0}\left( h \circ(\xi + 2i \overline{k}) \right) \circ(\xi + 2i \overline{k})  \ .
\end{eqnarray}
It will turn out that in the application to the Davey Stewartson II equation, the function $f$ is not needed, it just 
appears in the above derivation, but all the relevant information about the reflection coefficient is contained in the function $h$. 

It is useful to simplify the integral equation (\ref{eqh}) by first noting that
   \begin{eqnarray}
 &&   \mathcal{K}_{0}
 \left( h  \circ(\xi + 2i \overline{k}) \right) \circ(\xi + 2i \overline{k}) = 
  -i \mathcal{F} \left( 
    q \left[  \overline{ \mathcal{F}^{-1} \left( \frac{ h  }{\xi-2i \overline{k}}  \right)  }\right]  \right) \ , 
    \end{eqnarray}
so that the integral equation (\ref{eqh}) becomes 
\begin{eqnarray}
\label{eqh2}
h-\mathcal{K}_{0} \left( 
-i \mathcal{F} \left( 
    q \left[  \overline{ \mathcal{F}^{-1} \left( \frac{ h  }{\xi-2i \overline{k}}  \right)  }\right]  \right) 
\right) = - i \mathcal{F}(q) \ .
\end{eqnarray}
The above simplification follows from a straightforward calculation:
  \begin{eqnarray}
 &&   \mathcal{K}_{0}
 \left( h  \circ(\xi + 2i \overline{k}) \right) \circ(\xi + 2i \overline{k})
  = -i \mathcal{F} \left( 
    q \left[\overline{ \mathcal{F}^{-1} \left( \frac{ h  \circ(\xi + 2i \overline{k}) }{\xi}  \right) }\right]  \right)
    \circ(\xi +2 i \overline{k}) \\
    &&
     \ = \  -i \mathcal{F} \left( 
    q \left[\overline{ \mathcal{F}^{-1} \left( \frac{ h  }{\xi-2i \overline{k}}  \right)   \circ(\xi + 2i \overline{k})}\right]  \right) 
    \circ(\xi + 2i \overline{k})
    \\
    && = \  -i \mathcal{F} \left( 
    q \left[ e^{ kz - \overline{kz}} \overline{ \mathcal{F}^{-1} \left( \frac{ h  }{\xi-2i \overline{k}}  \right)  }\right]  \right) \circ(\xi + 2i \overline{k}) \\
    && = \ \ \ = \  -i \mathcal{F} \left( 
    q \left[  \overline{ \mathcal{F}^{-1} \left( \frac{ h  }{\xi-2i \overline{k}}  \right)  }\right]  \right) 
    \end{eqnarray}

\subsection{Algorithm to compute CGO solutions} \label{Subsec:5.3}
As in section \ref{Sec:3}, we introduce a standard discretisation of 
$x$, $y$ and the dual Fourier variables $\xi_{1}$, $\xi_{2}$. The 
role of the variable $k$ is in a sense also dual to $z$ as can be 
seen from equation (\ref{eq:MainEQa}). Thus we choose $2k$ on the 
same grid as $\xi$. The algorithm to compute the CGO solutions via 
(\ref{eqh2}) is thus as follows:
\begin{itemize}
    \item  For given $k$, approximate (\ref{eqh2}) as in section 
    \ref{Sec:3} by discretizing $z$ and $\xi$. 

    \item  Solve the resulting linear system for the resulting vector 
    $h$ with GMRES. 

    \item The task of computing the left hand side of (\ref{eqh2}) breaks into two parts:
    \begin{enumerate}
    \item Compute $\mathcal{F}^{-1} \left( \frac{ h  }{\xi-2i \overline{k}}  \right)  $ by the shifting and regularization argument explained in (\ref{eq:dbarinvShift}), and then compute 
    \begin{eqnarray}
    H:=
  -i \mathcal{F} \left( 
    q \left[  \overline{ \mathcal{F}^{-1} \left( \frac{ h  }{\xi-2i \overline{k}}  \right)  }\right]  \right) \ .
  \end{eqnarray}
  \item Finally, compute $\mathcal{K}_{0}(H)$ by using the previous regularization argument from Section 4.
   \end{enumerate}

\end{itemize}
The computational cost is roughly twice the one of the algorithm of 
section \ref{Sec:3} since the latter algorithm has to be applied 
twice (mainly 6 two-dimensional FFT per iteration). The advantage is 
that machine precision can be reached with much smaller values of 
$N_{x}$, $N_{y}$ for all values of $k$ with the present algorithm.

\section{Time dependence of DS II solutions}
In this section we study the time dependence of DS II solutions 
$q(x,y,t)$ for given initial data $q(x,y,0)$. 

{\bf Parameter Convention}:  As noted in Subsection \ref{Subsec:5.3}, it is convenient to choose $2k$ on the same grid as $\xi$.  In this Section, we make that choice, by replacing $k$ with $k/2$, so that $k$ and $\xi$ are dual variables and live on the same grid.

\subsection{Construction of solutions to the DS II equation via D-bar 
problems}

Recall that the D-bar equation associated to the DS II equation is the system (\ref{eq:021})-(\ref{eq:022}), normalized by (\ref{eq:022}).  Writing $S_{+}$ for $\xi \hat{M}_{1}$, and $S_{-}$ for $\xi \hat{M}_{2}$, we have two equations:
\begin{eqnarray}
\label{eq:MainEQDSpm}
S_{\pm}(\xi) =  \mp i \mathcal{F}\left\{  Q e^{ \frac{\overline{kz}-kz}{2}} \left(\overline{ \mathcal{F}^{-1} \left( \frac{1}{\xi} S_{\pm}(\xi) \right) }  \right) \right\} 
 \mp i  \mathcal{F}\left\{  Q e^{ \frac{\overline{kz}-kz}{2}}  \right\} \ .
\end{eqnarray}
For each equation, we seek to represent the solution in the form (\ref{eq:Sform}):
\begin{eqnarray*}
S_{\pm} = f_{\pm}(\xi) + h_{\pm}\circ(\xi + i \overline{k}) \ .
\end{eqnarray*}
It is straightforward to verify, following the same calculations from (\ref{eq:Sform}) to (\ref{eqh}), that
\begin{eqnarray}
\label{eq:hpmdef}
h_{\pm} = \pm h \ , 
\end{eqnarray}
where $h$ solves the equation
\begin{eqnarray}
h- \mathcal{K}_{0}\left(
\mathcal{K}_{0}\left( h  \circ(\xi + i \overline{k}) \right) \circ(\xi + i \overline{k}) 
\right) =- i \mathcal{F}\left( q \right) \ .
\end{eqnarray}
Moreover, 
\begin{eqnarray*}
f_{+} = f_{-}  = \mathcal{K}_{0}\left( h \circ(\xi + i \overline{k}) \right) \circ(\xi + i \overline{k})  \ .
\end{eqnarray*}

For a given value of $k\in \mathbb{C}$, the reflection coefficient is obtained from the computed function $h$ via
\begin{eqnarray}
\label{eq:rviah}
r(k) = i h(i \overline{k}) \ .
\end{eqnarray}
This can be seen by first writing a basic representation formula for $m_{\pm}$:
\begin{eqnarray}
m_{\pm} =  \frac{1}{ \pi} \iint_{\mathbb{C}} \frac{\dbar m_{\pm}}{z - z'}d^{2} z' \ , 
\end{eqnarray}
from which one may read off the $\displaystyle \mathcal{O} \left( \frac{1}{z} \right)$ term in the large $z$ expansion:
\begin{eqnarray}
\label{eq:m1def}
m_{\pm}^{(1)} =  \frac{1}{ \pi} \iint_{\mathbb{C}} \dbar m_{\pm} \ d^{2} z'  \ = \ \left.  \left( i \xi  \hat{m}_{\pm} \right) \right|_{\xi = 0} \ = \ i S_{\pm} (0) \ .
\end{eqnarray}
Now formula (\ref{eq:rviah}) for the reflection coefficient follows from (\ref{eq:m1def}) along with  (\ref{eq:Sform}), (\ref{eq:hpmdef}), and (\ref{eq:fpmdef}).  

\vskip 0.2in

To solve the DS II equation then, we compute 
as described above  the reflection coefficient 
(\ref{eq:rviah}) and obtain its time dependence from (\ref{rt}). It 
can be seen that the latter has a rather simple form which 
immediately implies that $|r|$ is constant in time which is 
equivalent to the existence of an infinite number of conserved 
quantities reflecting the complete integrability of the DS II 
equation. The most prominent conserved quantities are the $L^{2}$ 
norm of $q$ and the energy
\begin{equation}
    E=\iint_{\mathbb{R}^{2}}\left[|\partial_x q(x,y,t)|^2 - 
|\partial_y q(x,y,t)|^2    
+|q(x,y,t)|^{4}-\frac{1}{2}\left(\phi(x,y,t)^{2}+(\partial_{x}^{-1}\partial_{y}\phi(x,y,t))^{2}\right) \right]
      \label{energy}.
\end{equation}
It is well known that conserved quantities, the conservation of which 
is not implemented in the code, provide a test for numerical schemes. 
In  \cite{KR2011} it was shown that the error in the numerical 
conservation of the DS II energy tends to overestimate the numerical 
accuracy of the solution in an $L^{\infty }$ sense by two to three 
orders of magnitude.

It is a remarkable fact of integrable systems that the inverse 
scattering techniques, here a D-bar problem, immediately allow to go 
to the final time $t$ one is interested in. The standard approach is 
as follows: determine the reflection coefficient $r(k,0)$ as discussed in 
the previous section by solving (\ref{eqh}) for given $q(x,y,0)$. Then 
solve the same equation with $r(k,t)$ via (\ref{rt}) after 
interchanging $q$ and $r$, and $z$ and $k$. Thus in contrast to the 
direct numerical solution of DS II as in \cite{KR2011}, no 
intermediate time steps have to be computed. One just needs to solve 
twice integral equations of the form (\ref{eqh}) for all considered 
values of $k$ and $z$. 

\subsection{Reconstruction of the initial data at $t=0$}

A first test of this approach is to recover the potential $q(x,y,0)$ 
from the reflection coefficient $r(k,0)$. To this end one fixes the parameters 
$L_{x}$, $L_{y}$ in (\ref{domain}) and $N_{x}$, $N_{y}$ in 
(\ref{xidomain}). Maximal precision can only be reached for if both 
$q(x,y,0)$ and its Fourier transform decrease to machine precision at 
the boundaries of the considered domains. Thus if one wants to study 
the dependence of the difference between the computed $q(x,y,0)$ via 
the reflection coefficient and the original $q(x,y,0)$ on the 
parameters $L_{x}=L_{y}=L$ and $N_{x}=N_{y}=N$, the pairs $L$, $N$ 
have to be chosen such that both $q(x,y,0)$ and its Fourier transform 
decrease roughly to the same order of magnitude. If this is 
respected, then the numerical error 
decreases as expected exponentially, see the table in Fig.~\ref{table}. It shows 
that the initial data can be recovered to better than $10^{-13}$ with 
$N=256$ Fourier modes in each direction. 
\begin{figure}[htb!]
\begin{tabular}{|c|c|c|}
\cline{1-3}
Fourier Modes & L & Error\\
\hline
$2^{3}$ & 0.7515& \texttt{7.09e-03}  \\
\hline
$2^{4}$& 1.075& \texttt{3.1872e-04}\\
\hline
$2^{5}$ &1.5& \texttt{1.665e-06}\\
\hline
$2^{6}$& 2.1213& \texttt{1.736e-09}\\
\hline
$2^{7}$& 3.2& \texttt{2.40e-13}\\
\hline
$2^{8}$& 4.2& \texttt{5.0e-14}\\
\hline
\end{tabular}
\caption{Values of number $N$ of Fourier modes and $L$ for 
$M=11$ such that both the reconstructed potential $q(x,y,0)$
and the reflection coefficient decrease to roughly the same precision at the 
boundaries of the respective computational domains, and the 
corresponding errors in the reconstruction of the initial data 
$q(x,y,0)=\exp(-x^{2}-y^{2})$ after solving twice equations of the 
form (\ref{eqh}). }  
\label{table}
\end{figure}

\subsection{DS II solutions for $t>0$}

\begin{figure}[htb!]
   \includegraphics[width=0.49\textwidth]{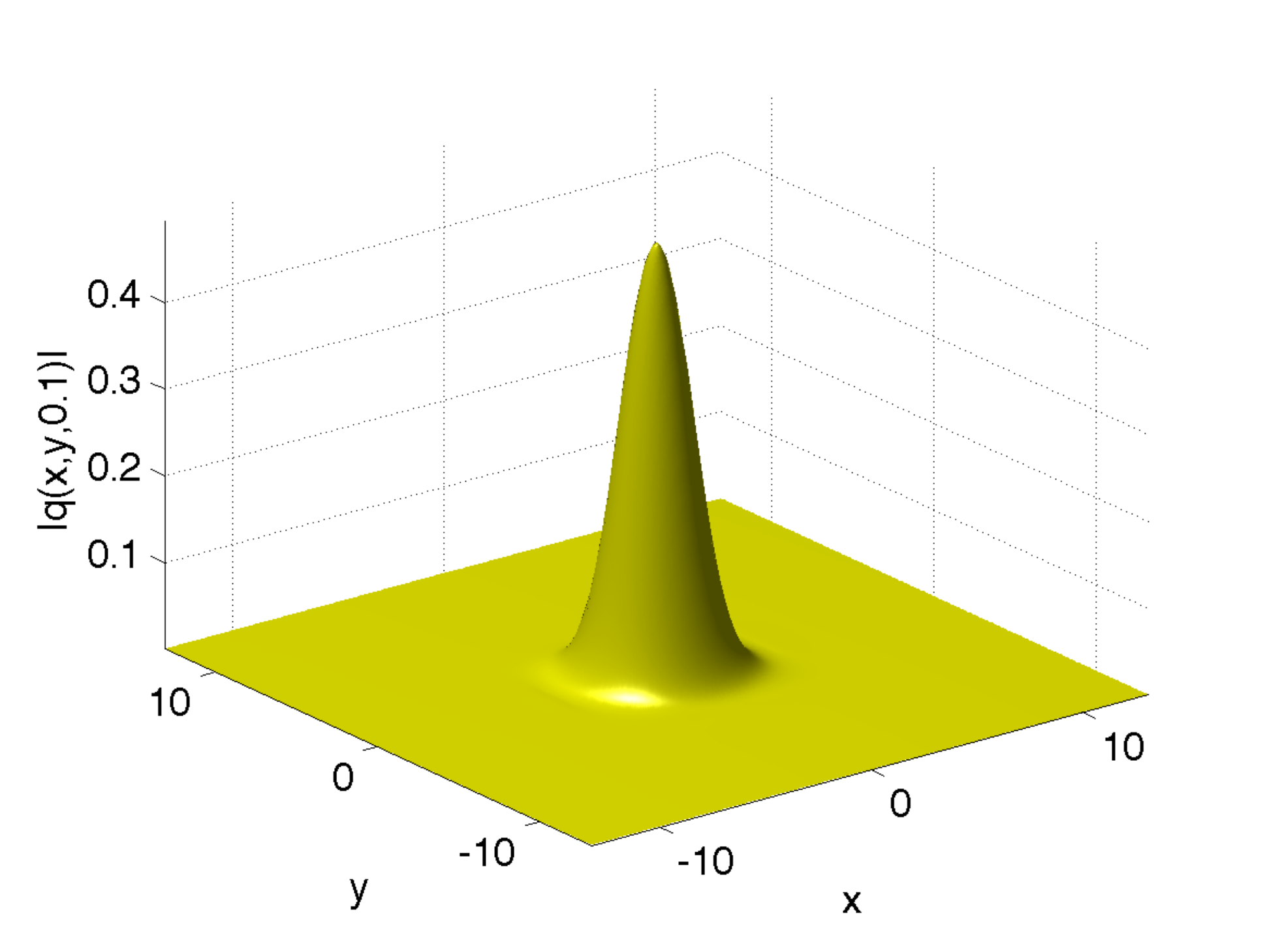}
   \includegraphics[width=0.49\textwidth]{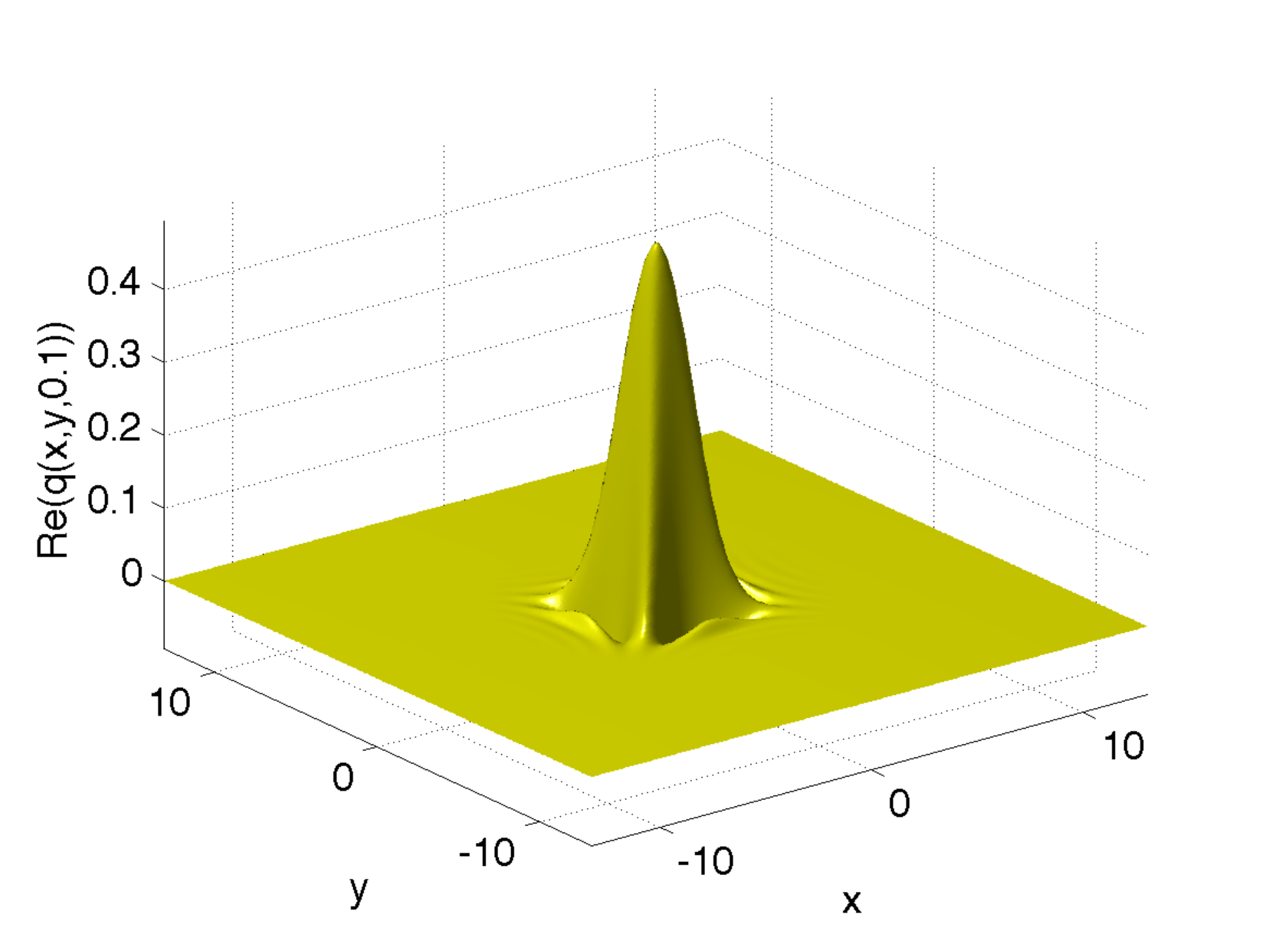}
\caption{Solution to the DS II equation for the initial data 
$q(x,y,0)=\exp(-x^{2}-y^{2})$ for $t=0.8$; on the left the modulus of 
the solution, on the right the real part. }  
\label{DSIIsol}
\end{figure}

It is known from the linear Schr\"odinger equation that Gaussian 
initial data are dispersed with time. A similar effect is to be 
expected in the case of the DS II solution. Thus parameters optimized 
for $t=0$ will in general not provide the same resolution for $t>0$. 
The reason is that in addition to the dispersive effects of DS 
II, there will be increasingly oscillations in the real and imaginary 
part of the DS II solution, see Fig.~\ref{DSIIsol}. Nonetheless, when we consider the final time $t=0.8$, 
the computed solution still conserves the energy to an error of order $10^{-14}$.
Of course, larger values of $t$ will require higher 
resolutions, which does not pose a principle obstacle, but might be 
best done on parallel computers, see the comments in the following 
section. 

The D-bar code can be used as a benchmark for direct solvers. For 
initial data in $\mathcal{S}(\mathbb{R}^{2})$, FFT techniques are 
again the appropriate method for such solvers. The DS II equation (\ref{eq:DSIIa}) and 
(\ref{eq:DSIIb}) can be written in Fourier space in the form
\begin{equation}
    i(\mathcal{F}q)_{t}-(\xi_{1}^{2}-\xi_{2}^{2})\mathcal{F}q-\mathcal{F}
    \left[\mathcal{F}^{-1}\left(\frac{\xi_{1}^{2}-\xi_{2}^{2}}{\xi_{1}^{2}+\xi_{2}^{2}}
    \mathcal{F}(|q|^{2})\right)q\right]=0
    \label{DSfourier};
\end{equation}
here equation (\ref{eq:DSIIb}) has been solved in Fourier space to 
eliminate the function $\varphi$ in (\ref{eq:DSIIa}). 
Approximating the Fourier transform as in the previous sections via a 
discrete Fourier transform, one ends up for (\ref{DSfourier}) with a 
finite dimensional system of ODEs in $t$ which are then numerically 
integrated with respect to $t$. 
The first approach along these lines appears to have been realised in 
\cite{WW}.  More recently, in \cite{KR2011}, fourth order stiff integrators for DS II 
have been studied. 

The inversion of the Laplace 
operator in  (\ref{eq:DSIIb}) has introduced the term 
$(\xi_{1}^{2}-\xi_{2}^{2})/(\xi_{1}^{2}+\xi_{2}^{2})$ in 
(\ref{DSfourier}). This term is much better behaved than the term 
$1/(\xi_{1}+i\xi_{2})$ appearing in the solution of the D-bar problem, 
but it is still a problem for a spectral approach. This can be easily 
seen by writing $\xi_{1}=|\xi|\cos \phi$ and $\xi_{2}=|\xi|\sin 
\phi$ which implies 
$(\xi_{1}^{2}-\xi_{2}^{2})/(\xi_{1}^{2}+\xi_{2}^{2})=\cos 2\phi$. 
Thus this factor does not have a well defined limit for 
$\xi_{1}=\xi_{2}=0$ and is put equal to zero there,  its average value with respect 
to $\phi$.
In the context of spectral accuracy as dicussed above, a non regular 
function  is a problem even if there is only one non-regular point 
since spectral methods are nonlocal. 

Due to the absence of an explicitly 
known smooth periodic DS II solution, spectral convergence had not 
been tested in \cite{KR2011}. The D-bar approach to DS II allows one to 
address this issue.  Using the direct solvers compared in 
\cite{KR2011}, with data and computational parameters consistent with the example of Fig.~\ref{DSIIsol}, using $N_{t}=10^{4}$ time steps, we find that the energy is conserved to the order of machine precision (the relative error between initial and final energy is $\texttt{2e-14}$). The difference between this 
numerical solution and the one obtained with the D-bar code is of the 
order of $10^{-6}$. This suggests that a regularization of the integrand in 
(\ref{DSfourier}) along the lines developed in this paper for the 
D-bar problem will be beneficial also in the context of direct 
numerical approaches for DS II. A related approach will be discussed elsewhere.

\begin{remark}
   One can also compare the performance of the D-bar solvers in 
   sections \ref{Sec:3} and \ref{Sec:4} in the context of DS II 
   solutions. Using the same parameters as 
   above for both codes, the following table shows the $L^{\infty}$ norm of the difference between both solutions.  Specifically, we numerically solve the inverse problem, with predetermined reflection coefficient $r_{0}(k_{1}, k_{2})e^{ -\frac{i t}{2} (\mbox{ Re}(k^{2}))}$, taking $t=0.1$, in two ways:  the only difference being in the way we computationally implement $\dbar^{-1}$, either using the method of Section \ref{Sec:3} or the method described in this section (using the ``shifting'' described in Section \ref{Sec:4}). 
   \begin{figure}[htb!]
\begin{tabular}{|c|c|c|}
\cline{1-3}
Fourier Modes & L & Error\\
\hline
$2^{6}$& 2.25& \texttt{7.3e-07}\\
\hline
$2^{7}$& 3.2& \texttt{3.4e-09}\\
\hline
$2^{8}$& 4.2& \texttt{5.2e-11}\\
\hline
\end{tabular}
\caption{$L^{\infty}$ norm of the difference between solution to DS-II at time $t=0.8$ obtained via the D-bar approach in two ways:  one way using the method of section \ref{Sec:3} for the inversion of the D-bar operator, and one way using the method of section \ref{Sec:4} for this inversion.  }  
\end{figure}
 This indicates that the more economic code 
   of section \ref{Sec:3} can be used in this context if machine 
   precision is not needed. The reason for the relatively high 
   accuracy of this code is that the biggest errors are introduced 
   there for $k$ and $z$ close to the boundaries of the computational 
   domains    where the reflection coefficient and the solution are 
   smallest. Thus these errors contribute comparatively little to the 
   final result. But if the CGO solution of the D-bar problem is needed 
   for values of $k$ close to the boundary with high precision, the 
   code of section \ref{Sec:4} has to be used.
\end{remark}

\section{Outlook}
We have presented in this paper algorithms for the numerical 
computation of  CGO solutions to  D-bar problems. The presented 
approach is based on FFT techniques and a regularization of the 
functions submitted to FFTs via $\phi$-functions. It was shown that the 
algorithms show spectral convergence and thus allow computations of 
the solution to the D-bar problem with essentially machine precision. 
The computational cost is mainly in two-dimensional FFTs which 
appear in an iterative resolution of the integral equation 
(\ref{eqh}) via GMRES. 

Whereas the convergence of the algorithm is as shown spectral, the 
solution of the equation (\ref{eqh}) has to be computed for each 
considered value of $k$ separately, which is time consuming. Thus a 
simple acceleration of the code on modern computers could be achieved 
by parallelizing the code in a way that the solution $h$ is computed 
at the same time for several values of $k$. If higher resolutions 
than $N_{x}N_{y}>2^{25}$ are needed, the corresponding FFTs are 
difficult to handle on a serial computer. Such resolutions can be 
necessary to study \emph{dispersive shocks}, i.e., highly oscillatory 
regions of the DS II solution appearing in the \emph{semiclassical} 
regime, see \cite{KR2013b}. In this case a parallelization also of 
the two-dimensional FFTs is necessary as in \cite{KR2013b}. The study 
of both types of parallelization will be the subject of further work. 

The algorithms studied in the present paper are optimized for rapidly 
decreasing smooth functions. For algebraically decreasing functions 
and functions with compact support being just piecewise continuous, 
different techniques are necessary to obtain spectral convergence. 
Possible approaches are multi-domain spectral methods as in \cite{BK} 
and references therein, best after introducing special coordinates as 
polar coordinates. This will be studied in an ensuing publication.

\section*{Acknowledgement}
This work was supported in part by the Marie Curie IRSES program 
RIMMP.  KM was supported in part by the National Science Foundation under grant DMS-1401268; he thanks the faculty and staff of the University of Burgundy for their hospitality during his stay as a visiting professor where part of this work has been completed, and for funding from the CNRS as a visiting researcher at the IMB.  KM also thanks Robert Indik for sage advice and a number of helpful discussions.

\end{document}